\documentclass[11pt]{amsart}
\usepackage{mathrsfs}
\usepackage{amssymb}
\usepackage[all,pdf]{xy}
\usepackage{titletoc}
\usepackage{amscd}
\usepackage{amsmath, amssymb}
\usepackage{amsfonts}
\usepackage[colorlinks,linkcolor=blue,citecolor=blue, pdfstartview=FitH]{hyperref}

\numberwithin{equation}{section}

\theoremstyle{plain}
\newtheorem{thm}{Theorem}[section]
\newtheorem{theorem}[thm]{Theorem}
\newtheorem{lemma}[thm]{Lemma}
\newtheorem{corollary}[thm]{Corollary}
\newtheorem{proposition}[thm]{Proposition}

\theoremstyle{definition}

\newtheorem{remark}[thm]{Remark}

\newtheorem{definition}[thm]{Definition}

\newtheorem{example}[thm]{Example}

\newtheorem{defn-thm}[thm]{Definition-Theorem}

\newcommand{\sA}{{\mathcal A}}

\newcommand{\sC}{{\mathcal C}}

\newcommand{\sF}{{\mathcal F}}
\newcommand{\sG}{{\mathcal G}}

\newcommand{\sJ}{{\mathcal J}}

\newcommand{\sL}{{\mathcal L}}
\newcommand{\sM}{{\mathcal M}}
\newcommand{\sN}{{\mathcal N}}
\newcommand{\sO}{{\mathcal O}}

\newcommand{\sS}{{\mathcal S}}
\newcommand{\sT}{{\mathcal T}}

\newcommand{\sX}{{\mathcal X}}
\newcommand{\sY}{{\mathcal Y}}

\newcommand{\B}{{\mathbb B}}
\newcommand{\C}{{\mathbb C}}
\newcommand{\D}{{\mathbb D}}

\newcommand{\F}{{\mathbb F}}

\newcommand{\K}{{\mathbb K}}
\renewcommand{\L}{{\mathbb L}}

\newcommand{\N}{{\mathbb N}}
\renewcommand{\P}{{\mathbb P}}
\newcommand{\Q}{{\mathbb Q}}
\newcommand{\R}{{\mathbb R}}
\newcommand{\T}{{\mathbb T}}

\newcommand{\X}{{\mathbb X}}

\newcommand{\Z}{{\mathbb Z}}
\newcommand{\id}{{  id}}

\newcommand{\Hom}{{\mathrm{Hom}}}
\newcommand{\Ker}{{\mathrm{Ker}}}
\newcommand{\Coker}{{\mathrm{Coker}}}

\newcommand{\btheorem}{\begin{theorem}}
\newcommand{\etheorem}{\end{theorem}}
\newcommand{\bproposition}{\begin{proposition}}
\newcommand{\eproposition}{\end{proposition}}
\newcommand{\bdefinition}{\begin{definition}}
\newcommand{\edefinition}{\end{definition}}
\newcommand{\bcorollary}{\begin{corollary}}
\newcommand{\ecorollary}{\end{corollary}}
\newcommand{\bproof}{\begin{proof}}
\newcommand{\eproof}{\end{proof}}
\newcommand{\bremark}{\begin{remark}}
\newcommand{\eremark}{\end{remark}}
\newcommand{\eexample}{\end{example}}
\newcommand{\bexample}{\begin{example}}

\newcommand{\elemma}{\end{lemma}}
\newcommand{\blemma}{\begin{lemma}}
  \newtheorem{thmx}{Theorem}

\newcommand{\ee}{\end{eqnarray*}}
\newcommand{\be}{\begin{eqnarray*}}

\newcommand{\beq}{\begin{equation}}
\newcommand{\eeq}{\end{equation}}

\newcommand{\bm}{\boldsymbol m}
\newcommand{\bn}{\boldsymbol n}
\newcommand{\bff}{\boldsymbol f}
\newcommand{\bg}{\boldsymbol g}

\newcommand{\bp}{\boldsymbol{p}}
\newcommand{\bee}{\boldsymbol{e}}
\newcommand{\bla}{\boldsymbol{\lambda}}

\newcommand{\lrw}{\longrightarrow}

\usepackage{fancyhdr}
\usepackage{tikz-cd}
\usepackage[all]{xy}
\begin{document}

\title{Quantum cluster algebras associated to weighted projective lines}
\author{Fan Xu, Fang Yang*}
\address{Department of Mathematical Sciences\\
Tsinghua University\\
Beijing 100084, P.~R.~China}
\email{fanxu@mail.tsinghua.edu.cn(F. Xu)}
\address{Department of Mathematical Sciences\\
Tsinghua University\\
Beijing 100084, P.~R.~China}
\email{yangfang19@mails.tsinghua.edu.cn(F. Yang)}

%
\thanks{$*$~Corresponding author.}


\begin{abstract}
  Let $\X_{\bp,\bla}$ be a weighted projective line. We define the quantum cluster algebra of $\X_{\bp,\bla}$ and realize its specialized version as the subquotient of the Hall algebra of $\X_{\bp,\bla}$ via the quantum cluster character map. Inspired by \cite{Chen2021}, we prove an analogue cluster multiplication formula between quantum cluster characters. As an application, we obtain the polynomial property of the cardinalities of Grassmannian varieties of exceptional coherent sheaves on $\X_{\bp,\bla}$ . In the end, we construct several bar-invariant $\Z[\nu^{\pm}]$-bases for the quantum cluster algebra of the projective line $\P^1$ and show how it coincides with the quantum cluster algebra of the Kronecker quiver.
\end{abstract}

\maketitle
\tableofcontents
\section{Introduction}
  The cluster algebras is a commutative algebra generated by a family of generators called cluster variables, which was introduced by Fomin and Zelevinsky \cite{SergeyFomin2002, Fomin2003} in order to study total positivity in algebraic groups and the specialization of canonical bases of quantum groups at $q = 1$. In \cite{Buan2006a}, Buan et al. introduces the cluster category as an additive categorification of the cluster algebra.   Cluster algebra and cluster category are closely related by the Caldero-Chapoton map in \cite{Caldero2006a} and the Caldero-Keller multiplication theorem in \cite{Caldero2008,Caldero2006a}. Caldero and Keller \cite{Caldero2008} proved  the following formula (called $cluster$ $multiplication$ $formula$) 
   \begin{equation}\label{eqn1.1}
      \chi(\P \operatorname{Ext}^1(M,N))X_MX_N=\sum_E (\chi(\P \operatorname{Ext}^1(M,N)_E)+\chi(\P \operatorname{Ext}^1(N,M)_E))X_E.
   \end{equation} 
  for any objects $M, N \in \sC_Q$ such that $\operatorname{Ext}^1_{\sC_Q}(M, N)\neq 0$ for $Q$ is of finite type.  And Caldero-Keller \cite{Caldero2006a} showed 
     \begin{equation}\label{eqn1.2}
      X_MX_N=X_E+X_{E'}.
     \end{equation}
  for $M,N\in \sC_Q$ indecomposable such that  $\operatorname{Ext}^1_{\sC_Q}(M,N)$ is one-dimensional. Various generalizations of the above formulas were made by Hubery \cite{Hubery2010}, by Xiao and Xu \cite{Xiao2010,Xu2010}, by Fu and Keller \cite{Fu2010} and by Palu \cite{Palu2008,Palu2012}. In the cluster theory, the Caldero-Chapoton map and the cluster multiplication theorem play a very important role in proving some structural results such as bases with good properties, positivity conjecture, denominator conjecture and so on (cf. \cite{Caldero2008,Ding2013}).

  As a quantum analogue of cluster algebras, quantum cluster algebras were defined by Berenstein and Zelevinsky \cite{Berenstein2005} in order to study canonical bases for quantum groups of Kac-Moody type. Under the specialization $q = 1$, the quantum cluster algebras are exactly cluster algebras. As for the quantum cluster algebra of a valued acyclic quiver, Rupel \cite{Rupel2011} defined a quantum analogue of the Caldero-Chapoton map over a finite field. The quantum version of Equation (\ref{eqn1.2}) was proved by Rupel in \cite{Rupel2011} for indecomposable rigid objects for all finite type valued quivers, by Qin \cite{Qin2012} for indecomposable rigid objects for acyclic quivers.  Chen-Ding-Zhang \cite{Chen2021} gave the cluster multiplication formulas between any two quantum cluster characters. These formulas were a quantum version of the cluster multiplication formula in Equations (\ref{eqn1.1})  and (\ref{eqn1.2}) for acyclic quantum cluster algebras. 

  In \cite{Ding2020a}, Ding-Xu-Zhang realized an acyclic quantum cluster algebra as a subquotient of certain derived Hall algebra. This result was refined and generalized by Fu-Peng-Zhang \cite{Fu2020}  via the integration map from the Hall algebra of an acyclic quiver to certain quantum torus. This provides a connection between Hall algebras and quantum cluster algebras. Then one may define a ``new" quantum cluster algebra as a proper subquotient of the Hall algebra. As shown by Kapranov \cite{Kapranov1997} and then Schiffmann \cite{Schiffmann2004}, the Hall algebra of a weighted projective line gives a categorification of the positive part of the associated quantum loop algebra. If we can define a kind of quantum cluster algebra as a subquotient of the Hall algebra of weighted projective lines, then it may be possible to study the canonical bases of quantum loop algebras by using the quantum cluster algebras of weighted projective lines. 

  The aim of this paper is to define the quantum cluster algebra associated to a weighted projective line $\X_{\bp,\bla}$ by the Hall algebra of the category $\mathrm{Coh}(\X_{\bp,\bla})$ of coherent sheaves on $\X_{\bp,\bla}$. To define a cluster algebra as a subalgebra of the quantum torus generated by some elements indexed by a set $\sJ$ of indecomposable rigid objects  in certain cluster category $\sC$ over an algebraically closed field, it is required that $\sJ$ admits a cluster structure (see \cite[Section 1]{Buan2010}). But for the definition of a quantum cluster algebra, the first difficulty is to find a cluster structure independent of finite fields. In this paper, we use the valued regular $m$-tree $\T_m(k)$ (defined in Definition \ref{def3.6}) to denote the cluster structure of the cluster category $\sC(\mathrm{Coh}(\X_{\tilde{\bp},\tilde{\bla}})_k)$. Besides, we also need to find a suitable compatible pair $(\varLambda, \tilde{B})$ (see \cite[Section 3]{Berenstein2005}). 
  
  Our strategy is firstly to show there is a common valued regular tree $\T_m$ over finite fields $\F_{q^r}$ and the algebraic closure $\bar{\F}_q$ for some fixed prime $q$, whose proof will be given in Appendix \ref{ap}. Then to show that the valued regular trees over algebraic closures of distinct finite characteristics are the same by taking use of quiver with potentials. For the skew-symmetrizable matrix $\tilde{B}$, as in the case of acyclic quivers, we let $\tilde{B}$ be the skew-symmetric Euler form on $\tilde\sA:=\mathrm{Coh}(\X_{\tilde{\bp},\tilde{\bla}})$, where each item of $\tilde{\bp}$ is odd. Due to \cite[Section 9]{Geigle1991}, the category $\mathrm{Coh}(\X_{\bp,\bla})$ can be embedded into $\mathrm{Coh}(\X_{\tilde{\bp},\tilde{\bla}})$ if $\bp\leq \tilde{\bp}$ and $\lambda=\tilde{\lambda}$, which makes sure  that the principal part of $\tilde{B}$ is the skew-symmetric Euler form on $\mathrm{Coh}(\X_{\bp,\bla})$. The definition of quantum cluster algebra $\sA(\varLambda,\tilde{B}(\bp,\bla))$ of $\X_{\bp,\bla}$ is given in Section \ref{sec3.4}.  The first difference between the quantum cluster algebras of weighted projective line and the one of acyclic quiver is that the exchange matrix of the $1$-th mutation from the initial cluster-tilting object generally may not be determined by the initial compatible pairs of the quantum cluster algebras of $\X_{\bp,\bla}$. The essential reason is that the cluster category of a weighted projective line $\X_{\bp,\bla}$ may not be triangle equivalent to the cluster category of an acyclic quiver except for domestic type (see \cite[Remark 5.4]{Geigle1987}). Hence wo do not know whether the quantum cluster algebras of weighted projective lines admit Laurent phenomenon in general. The second difference between them is that the skew-symmetric form $\varLambda$ of the quantum cluster algebra of $\X_{\bp,\bla}$  does not change after mutations.
  
  In Section \ref{sec3}, we construct an algebra homomorphism (called quantum cluster character map) $X_?$ from the $\varLambda$-twisted Hall algebra $H_{\varLambda}(\tilde{\sA})$ to the specialized complete quantum torus $\hat{\sT}_{\varLambda,v}$, then in Section \ref{sec4.1} show a quantum analogue of the cluster multiplication formula (\ref{eqn1.1}) in $\hat{\sT}_{\varLambda,v}$ as Chen-Ding-Zhang did in \cite{Chen2021}:
  \begin{thmx}[Theorem 
  \ref{thm3.2}]
      For $\sM, \sN\in \tilde{\sA}$, in  $\hat{\sT}_{\varLambda,v}$ we have:
               \begin{align*}
                  &(q^{[\sM,\sN]^1}-1)X_{\sM}X_{\sN}=q^{\frac{1}{2}\varLambda(\bm^*,\bn^*)}\sum_{\sL\neq [\sM\oplus \sN]} |\operatorname{Ext}^1_{\tilde{\sA}}(\sM,\sN)_{\sL}| X_{\sL}\\
                  &\ \ \ \ \ \ +\sum_{[\sG],[\sF]\neq [\sN]} q^{\frac{1}{2}\varLambda((\bm-\bg)^*,(\bn+\bg)^*)+\frac{1}{2} \langle \bm-g,\bn \rangle }|_{\sF}\Hom_{\tilde{\sA}}(\sN,\tau\sM)_{\tau\sG}| X_{\sG}X_{\sF},
               \end{align*}
    \end{thmx}
  
  As an application, it is proved in Section \ref{sec4.2} that for an indecomposable rigid object $T_i(t)$ for $t\in \T_m$, there is a $\Z$-polynomial $P(z)$ such that the cardinality of $\mathrm{Gr}_{\boldsymbol{e}}(T_i(t)^k)$ is $P(|k|^{\frac{1}{2}})$. As a result, the generators of  the quantum cluster algebra $\sA(\varLambda,\tilde{B}(\bp,\bla))$ defined recursively by mutation formulas can be described as certain quantum cluster characters as stated in the following
  \begin{thmx}[Theorem \ref{thm4.6}]
    The quantum cluster algebra $\sA(\varLambda,\tilde{B}(\bp,\bla))$ as a subalgebra of $\hat{\sT}_{\varLambda}$ is generated by $X_{T_i(t)}$ for $t\in \T_n(\bp,\bla)$ and $X_{T_l(t_0)}^{\pm}$ for $n< l\leq m$. 
  \end{thmx}
   
  As another application of the cluster multiplication formula, in Section \ref{sec4.3} we show that the specialized quantum cluster algebra $\sA_q(\varLambda,B(\tilde{\bp},\tilde{\bla}))$ is a subquotient of the $\varLambda$-twisted Hall algebra $H_{\varLambda}(\tilde{\sA}_k)$. We prove the following
  \begin{thmx}[Theorem \ref{thm4.9}]
  There is an isomorphism of algebras :
  $$\phi_k: (\mathrm{CH}'_{\varLambda}(\tilde{\sA}_k)\otimes_{\Z[v^{\pm1}]} \Q)/I_k  \lrw \sA_q(\varLambda,B(\tilde{\bp},\tilde{\bla}))\otimes_{\Z[v^{\pm1}]} \Q,$$
  which maps $[T_i(t)^k]$ to $X_{T_i(t)^k}$ for $1\leq i\leq m$ and $t\in \T_m$.
  \end{thmx}

   In Section \ref{sec5}, we study the quantum cluster algebra $\sA(\varLambda,B)$ of the projective line $\P^1$ and show how it coincides with the quantum cluster algebra of the Kronecker quiver. We obtain
   \begin{thmx}[Theorem \ref{thm5.14}]
     Each one of the following sets gives rise to a bar-invariant $\Z[\nu^{\pm1}]$-basis for $\sA(\varLambda,B)$:
     \begin{center}
       $ \B_1^{tor}\cup \bar{\B}^{vet}$, $\ \ $  $\B_2^{tor}\cup \bar{\B}^{vet}$, $\ \ $   $\B_3^{tor}\cup \bar{\B}^{vet}$.
     \end{center}
    \end{thmx}
    These bases above are corresponding to the bar-invariant $\Z[\nu^{\pm1}]$-bases of the quantum cluster algebra $\sA(2,2)$ of Kronecker quiver constructed by Ding-Xu \cite{Ding2012}. They showed that under the specialization $\nu = 1$, these $\Z[\nu^{\pm1}]$-bases are exactly the canonical basis, semicanonical basis and dual semicanonical basis of the corresponding cluster algebra. 

 \subsection*{Acknowledgement} The authors are supported by the NSF of China (No.12031007) The authors thank Prof. Changjian Fu for pointing out that the cluster category of weighted projective line has non-degenerate quiver with potential. The authors also appreciate Prof. Xueqing Chen for helpful comments. The second author also thanks Prof. Yu Zhou for kindly answering many questions about the tilting theory. 

\subsection*{Conventions}
  Throughout this paper, denote by $k$ a finite field. $K$ is denoted to be an algebraically closed field of finite characteristic.  Let $\nu$ be a formal variable.  Denote by $\D=\Hom_k(-,k)$ the $k$-duality. $\sA_k$ (resp. $\tilde{\sA}_k$) is the category of coherent sheaves on weighted projective line $\X_{\bp,\bla}$ (resp. $\X_{\tilde{\bp},\tilde{\bla}}$) over $k$. We will omit $\sA_k$ for $\sA$ when it does not cause any confusions. In the Hall algebra of $\sA$, denote by $[\sF]$ the isoclass of $\sF$. In the Grothendieck group $K_0(\sA)$, we will denote by $\hat{\sF}$ (some times by $[\sF]$) the class of $\sF\in \sA$. Let $\sC(\sA)$ be the cluster category of $\sA$. every cluster-tilting object is assumed to be basic.  Let $\{\bee_i|\ 1\leq i\leq m\}$ be the canonical basis for $\Z^m$.
  \begin{itemize}
    \item $\sT_{\varLambda}$: the quantum torus associated to $\varLambda$,
    \item $\hat{\sT}_{\varLambda}$: the complete quantum torus associated to $\varLambda$,
    \item $\hat{\sT}_{\varLambda,v}$: the complete quantum torus specialized at $\nu=v$.
    \item $\sA(\varLambda,\tilde{B}(\bp,\bla))$: the quantum cluster algebra of $\X_{\bp,\bla}$,
    \item $\sA_q(\varLambda,\tilde{B}(\bp,\bla))$: the quantum cluster algebra of $\X_{\bp,\bla}$ specialized at $\nu=q^{\frac{1}{2}}$.
  \end{itemize}
 
\maketitle
\section{Preliminary}
 \subsection{Weighted projective lines}
   Let $k$ be a finite field $\F_q$ with $|k|=q$. Set $\boldsymbol{p}=(p_1,\cdots,p_N)$ be a collection of $N\geq 3$ positive integers. Denote by $S(\boldsymbol{p})$ the polynomial ring $k[X_1,\cdots,X_N]$ and consider the ideal $I(\bp,\boldsymbol{\lambda})$ generated by $X_i^{p_i}=X_2^{p_2}-\lambda_{i}X_1^{p_1}$ for $i\geq 3$, where $\lambda_1,\lambda_2,\cdots,\lambda_N$ are distinct points of $\P^1$ normalized in such a way that  $\lambda_1=\infty$, $\lambda_2=0$ and $\lambda_3=1$. Let $S(\boldsymbol{p},\boldsymbol{\lambda})$ be the quotient $S(\boldsymbol{p})/I(\bp,\boldsymbol{\bla})$.  Then $S(\boldsymbol{p},\boldsymbol{\lambda})$ is naturally graded by an abelian group $L(\boldsymbol{p}):=\Z\vec{x}_1\oplus \Z\vec{x}_2\cdots \oplus \Z\vec{x}_N/ (p_i\vec{x}_i-p_j\vec{x}_j,\forall i,j)$, and $X_i$ is associated with degree $\vec{x}_i$. Note that $S(\boldsymbol{p})$ is $L(\boldsymbol{p})$-graded and $I(\bp,\boldsymbol{\lambda})$ is generated by homogeneous elements, hence $S(\boldsymbol{p},\boldsymbol{\lambda})$ is also $L(\boldsymbol{p})$-graded. Denote $\vec{c}\in L(\bp)$ by $p_i\vec{x}_i$. The weighted projective line $\X_{\bp,\bla}$ is defined to be the spectrum $\mathrm{Spec}_{L(\bp)}S(\bp,\bla)$.

   Let $\mathrm{Coh}(X_{\bp,\bla})$ be the category of coherent sheaves on the weighted projective line $\X_{\bp,\bla}$, which is an abelian and hereditary category admitting an automorphism 
       $$\tau:\mathrm{Coh}(X_{\bp,\bla})\to \mathrm{Coh}(X_{\bp,\bla}),\ \sF\mapsto \sF(\vec{w}).$$
   where $\vec{w}=(N-2)\vec{c}-\sum_{i=1}^N \vec{x}_i=-2\vec{c}+\sum_{i=1}^{N} (p_i-1)\vec{x}_i\in L(\bp)$. Let $\mathrm{Vec}(X_{\bp,\bla})$ be the subcategory of $\mathrm{Coh}(X_{\bp,\bla})$ of locally free sheaves, and $\mathrm{Tor}(X_{\bp,\bla})$ be the subcategory of torsion sheaves. Since every coherent sheaf can be decomposed into a direct sum of a torsion part and a locally free part, we have
                  $$\mathrm{Coh}(X_{\bp,\bla})=\mathrm{Vec}(X_{\bp,\bla})\oplus \mathrm{Tor}(X_{\bp,\bla})$$

   Set $\Lambda:=\{\lambda_1,\cdots,\lambda_N\}$. these points are called exceptional points on $\P^1$. For any $x\in\P^1$, let $\mathrm{Tor}_x$ be the subcategory of torsion sheaves supported on $x$. 

   \begin{lemma}[\cite{Schiffmann2012a}]\label{lem2.1}
    The category $\mathrm{Tor}(X_{\bp,\bla})$ decomposes as a direct product of orthogonal blocks
                    $$\mathrm{Tor}(X_{\bp,\bla })=\prod_{x\in \P^1-\Lambda} \mathrm{Tor}_x\times \prod_{i=1}^N \mathrm{Tor}_{\lambda_i}.$$
    Moreover, $\mathrm{Tor}_x$ is equivalent to the category $\mathrm{Rep}^{nil}_{k_x}A_0^{(1)}$ of nilpotent representations of the Jordan quiver over the residue field $k_x$, and $\mathrm{Tor}_{\lambda_i}$ is equivalent to the category $\mathrm{Rep}^{nil}_k A_{p_i-1}^{(1)}$ of nilpotent representations of the cyclic quiver $A_{p_i-1}^{(1)}$ over $k$.
   \end{lemma}

   Hence, we denote by $S_x\in \mathrm{Tor}_x$ the simple torsion sheaf corresponding to the simple module of $\mathrm{Rep}^{nil}_{k_x}A_0^{(1)}$ for $x\in \P^1-\Lambda$,  and by $S_{ij}$ the simple torsion sheaf corresponding to the simple module on the $j$-th vertex of $\mathrm{Rep}^{nil}_k A_{p_i-1}^{(1)}$, $1\leq i\leq N$, $1\leq j\leq p_i$. Denote  by $\sO$ the structure sheaf on $\X_{\bp,\bla}$.
   
   \begin{lemma}[\cite{Schiffmann2012a}]\label{lem1}
    The Grothendieck group $K_0(X_{\bp,\bla})$ of $\mathrm{Coh}(X_{\bp,\bla})$ is isomorphic to 
                    $$(\Z[\sO]\oplus \Z[S_x]\oplus \bigoplus_{1\leq i\leq N,1\leq j\leq p_i} \Z[S_{ij}])\big / J.$$
    where $J$ is the subgroup generated by $[S_x]-\sum_{j=1}^{p_i} [S_{i,j}]$ for $1\leq i\leq N$.
   \end{lemma}
    
   As a corollary, we have that
                 $$K_0(X_{\bp,\bla})\cong \Z[\sO]\oplus \Z[S_x]\oplus \bigoplus_{1\leq i\leq N,2\leq j\leq {p_i}} \Z[S_{ij}].$$

\subsection{ The Hall algebra of  \texorpdfstring{$\mathrm{Coh}\X_{\bp,\bla}$}{Lg}}
   Fix $\bp=(p_1,\cdots,p_N)$ and $\bla=(\lambda_1,\cdots,\lambda_N)$, we get a weighted projective line $\X_{\bp,\bla}$. Let $k=\F_q$. Denote by $\sA$ the category $\mathrm{Coh}(X_{\bp,\bla})_k$ over $k$ and $\mathrm{Iso}(\sA)$ the set of isoclasses of objects in $\sA$. Let $ \langle , \rangle $ be the Euler form of $\sA$ on the Grothendieck group $K_0(\sA)$, that is,
                          $$ \langle \hat{\sF},\hat{\sG} \rangle =\dim_k \Hom_{\sA}(\sF,\sG)-\dim_k \operatorname{Ext}^1_{\sA}(\sF,\sG).$$
   where $\sF,\sG\in \sA$ and $\hat{\sF}\in  K_0(\sA)$ represents the class of $\sF$.  The symmetric Euler form is given by $(\hat{\sF},\hat{\sG}):= \langle \hat{\sF},\hat{\sG} \rangle + \langle \hat{\sG},\hat{\sF} \rangle $.
             
  To simplify notations, we will write $[\sF,\sG]^0$ for $\dim_k \Hom_{\sA}(\sF,\sG)$ and $[\sF,\sG]^1$ for $\dim_k\operatorname{Ext}^1_{\sA}(\sF,\sG)$. Denote $g_{\sF,\sG}^{\sL}=\#\{\sL_1\subset \sL| \sL_1\cong \sG, \sL/\sL_1\cong \sF\}$.
             
  The dual Hall algebra $H^{\vee}(\sA)$ of $\sA$ is defined to be the $\Q$-vector space $\bigoplus\limits_{[\sF]\in \mathrm{Iso}(\sA)} \Q[\sF]$ equipped with the multiplication
          $$[\sF][\sG]:=\sum_{[\sL]}q^{ \langle \sF,\sG \rangle } \frac{|\operatorname{Ext}_{\sA}^1(\sF,\sG)_{\sL}|}{|\Hom_{\sA}(\sF,\sG)|} [\sL].$$
  In the sequel, we will write $f_{\sL}^{\sF,\sG}$ for $ \frac{|\operatorname{Ext}_{\sA}^1(\sF,\sG)_{\sL}|}{|\Hom_{\sA}(\sF,\sG)|}$.
  \bremark\label{rem1}
   (i)Note that we define the Hall algebra $H^{\vee}(\sA)$ by using another multiplication which is dual to the usual Hall multiplication by counting subobjects.  The original Hall algebra $H(\sA)$ is the $\Q$-vector space  $\bigoplus\limits_{[[\sF]]\in \mathrm{Iso}(\sA)} \Q[[\sF]]$ equipped with the multiplication
         $$[[\sF]][[\sG]]:=\sum_{[\sL]}q^{ \langle \sF,\sG \rangle } g_{\sF,\sG}^{\sL} [[\sL]].$$
             
   (ii) The category $\sA$ of coherent sheaves does not satisfy the finite subobject condition. For example,  the structure sheaf $\sO$ has subobjects $\sO(r\vec{c})$ for $r<0$. Hence, if we want to give a  comultiplication $\Delta:H^{\vee}(\sA)\to H^{\vee}(\sA)\hat{\otimes} H^{\vee}(\sA)$, then $H^{\vee}(\sA)\hat{\otimes} H^{\vee}(\sA)$ is simply the space of all formal (may be infinitely many) linear combinations of $[\sF]\otimes [\sG]$.
   \eremark
             
 \begin{lemma}[\cite{Schiffmann2012a}]
  The following defines on $H^{\vee}(\sA)$ the structure of a topological coassociative coproduct:
         $$\Delta([\sF])=\sum_{\sF_1,\sF_2}q^{ \langle \sF_1,\sF_2 \rangle } g_{\sF_1,\sF_2}^{\sF} [\sF_1]\otimes [\sF_2].$$
   \end{lemma}
             
 Define the twisted multiplication on $H^{\vee}(\sA)\hat{\otimes} H^{\vee}(\sA)$  by 
        $$([\sF_1]\otimes [\sF_2])([\sG_1]\otimes [\sG_2]):=q^{(\sF_2,\sG_1)+ \langle \sF_1,\sG_2 \rangle } [\sF_1][\sG_1]\otimes [\sF_2][\sG_2].$$
\begin{lemma}[\cite{Schiffmann2012a}]
  The comultiplication $\Delta: H^{\vee}(\sA)\lrw H^{\vee}(\sA)\hat{\otimes} H^{\vee}(\sA)$ is a homomorphism of algebras.
\bproof
    We have
      $$\Delta([\sF][\sG])=\sum_{\sL}q^{ \langle \sF,\sG \rangle }f^{\sF,\sG}_{\sL}\sum_{\sL_1,\sL_2} q^{ \langle \sL_1,\sL_2 \rangle }g_{\sL_1,\sL_2}^{\sL}[\sL_1]\otimes [\sL_2].$$
  On the other hand,
      $$\Delta([\sF])\Delta([\sG])=\sum_{\sL_1,\sL_2}\sum_{\sF_i,\sG_i}q^ag_{\sF_1,\sF_2}^{\sF}g_{\sG_1,\sG_2}^{\sG} f_{\sL_1}^{\sF_1,\sG_1}f_{\sL_2}^{\sF_2,\sG_2}[\sL_1]\otimes [\sL_2].$$
  where $a= \langle \sF_1,\sF_2 \rangle + \langle \sG_1,\sG_2 \rangle +(\sF_2,\sG_1)+ \langle \sF_1,\sG_2 \rangle + \langle \sF_1,\sG_1 \rangle + \langle \sF_2,\sG_2 \rangle = \langle \sF,\sG \rangle + \langle \sL_1,\sL_2 \rangle $ $- \langle \sF_1,\sG_2 \rangle $. To show $\Delta([\sF][\sG])=\Delta([\sF])([\sG])$, it suffices to show that 
       $$\sum_{\sL}f^{\sF,\sG}_{\sL}g_{\sL_1,\sL_2}^{\sL}=\sum_{\sF_i,\sG_i} q^{- \langle \sF_1,\sG_2 \rangle }g_{\sF_1,\sF_2}^{\sF}g_{\sG_1,\sG_2}^{\sG} f_{\sL_1}^{\sF_1,\sG_1}f_{\sL_2}^{\sF_2,\sG_2}.$$ 
  for any $\sL_1$ and $\sL_2$, which is precisely the  Green's formula in the \cite[Theorem 2]{Green1995}.
  \eproof
  \end{lemma}

 \section{Quantum cluster characters}\label{sec3}
   \subsection{Compatible pairs}\label{sec3.1}
     For the weighted projective line $\X_{\bp,\bla}$. Recall that the Grothendieck group $K_0(\sA)$ is isomorphic to $\Z[\sO]\oplus \Z[S_x]\oplus \bigoplus_{i,2\leq j\leq {p_i}} \Z[S_{ij}]\cong\Z^n$  by the corollary of Lemma \ref{lem1}, here $n=2+\sum_{i=1}^N (p_i-1)$. Note that for any element $\hat{\sF}$ in $K_0(\sA)$, we will write $\underline{\dim}\sF$ for the dimension vector of $\hat{\sF}$ under the basis $\mathbf{b}(p,\lambda)$:
         $$\{\hat{\sO},\hat{S_x},\hat{S}_{i,j} |1\leq i\leq N, 2\leq j\leq p_i\}$$
    Namely,
     $$\underline{\dim}\sF=a_1\hat{\sO}+a_2\hat{S_x}+a_3\hat{S}_{1,2}+\cdots+a_{p_1+1}\hat{S}_{1,p_1}+a_{p_1+2}\hat{S}_{2,2}+\cdots+a_n\hat{S}_{N,p_N}.$$

     Let $E:=E(\bp,\bla)$ be the $n\times n$ matrix associated to the Euler bilinear form $ \langle , \rangle $ such that 
            $$(\underline{\dim}\sF)^tE \underline{\dim}\sG= \langle \hat{\sF},\hat{\sG} \rangle .$$
    Denote $E^t$ by the transpose of E. Set $B(\bp,\bla):=E^t-E$. Then by direct computation, $ \langle \hat{\sO},\hat{S_{ij}} \rangle =\delta_{p_i,j}$ and $ \langle \hat{S_{ij}},\hat{\sO} \rangle =-\delta_{1,j}$ for $1\leq j\leq p_i$, the matrix $B(\bp,\bla)$ has the form: 
          $$\begin{bmatrix}
            B_0       &C_1 &C_2 &\cdots &C_N\\
            -C_1^t  &B_1 &0     &\cdots &0\\
            -C_2^t  &0    &B_2  &\cdots &0\\
            \vdots   &0   &0      &\ddots &0\\
            -C^t_N  &0   &0      &\cdots &B_N\\
          \end{bmatrix}$$
      where $B_0=\begin{bmatrix} 0 &-2\\2 &0     \end{bmatrix}$, matrix $C_i=\begin{bmatrix} 0 &0 &\cdots &-1\\ 0 &0 &\cdots &0\end{bmatrix}$, and $B_i$ is a square matrix of $p_i-1$ as follows:
          $$\begin{bmatrix}
           0 &1 &0 &\cdots  &0  &0 \\
           -1 &0 &1 &\cdots &0&0 \\
            0  &-1  &0 &\cdots &0 &0 \\
            \vdots & & &\ddots &  &\\
              0 &0 &0 &\cdots  &0  &1 \\
             0 &0 &0 &\cdots &-1&0 \\
          \end{bmatrix}$$
    Since $\det{B}(\bp,\bla)$ is the product of $\det{B_i}$, and $B_i$ is invertible iff $p_i-1$ is even for $i=1,\cdots, N$, then $B(\bp,\bla)$ is invertible if and only if all $p_i$ is odd. If $B$ is not invertible, we can embed $B$ into some $m\times m$ invertible matrix $\tilde{B}(\bp,\bla)=\tilde{E}^t-\tilde{E}$ such that $E$  is the upper submatrix of $\tilde{E}$.  
     
    In the following, we give the construction of $\tilde{B}(\bp,\bla)$ such that it is the matrix of skew-symmetric Euler form of another weighted projective line up to a choice of basis for its Grothendieck group. Without loss of generality, we assume that only $B_1$ are noninvertible. Hence $p_1$ is even. Set $\tilde{\bp}=(p_1+1,p_2,\cdots,p_N)$ and $\tilde{\bla}=\bla$. By \cite[Theorem 9.5]{Geigle1991}, if $\tilde{\bp}=(p_1+1,p_2,\cdots,p_N)$, then there exists an exact equivalence 
    $$\phi_{*}:\tilde{\sA}\big/{{\mathrm{add}\tilde{S}_{1,1}}}\simeq {\sA},$$
   such that $\phi_{*}(\tilde{\sO})=\phi_{*}(\tilde{\sO}(\vec{x}_1))=\sO$, $\phi_{*}(\tilde{S}_{i,j})=S_{i,j}$ if $i\neq 1$ and $\phi_{*}(\tilde{S}_{1,j})=S_{1,j-1}$ for $2\leq j\leq p_1+1$.  By direct computations, we have that $E(\bp,\bla)$ is the upper-left submatrix of $E(\tilde{\bp},\tilde{\bla})^*$, where $E(\tilde{\bp},\tilde{\bla})^*$ is the matrix of Euler form on $K_0(\tilde{\sA})$ under the following basis $\mathbf{b}(\tilde{\bp},\tilde{\bla})^*$:
    $$[\tilde{\sO}],[\tilde{S}_x], [\tilde{S}_{1,3}],\cdots,[\tilde{S}_{1,{p_1+1}}], [\tilde{S}_{2,2}],[\tilde{S}_{2,3}],\cdots, [\tilde{S}_{N,p_N}], [\tilde{S}_{1,2}].$$
    
    Set $\tilde{B}(\bp,\bla)$ to be $E(\tilde{\bp},\tilde{\bla})^{*t}-E(\tilde{\bp},\tilde{\bla})^*$, which is obtained from $B(\tilde{\bp},\tilde{\bla})$ by base change from $\mathbf{b}(\tilde{\bp},\tilde{\bla})$ to  $\mathbf{b}(\tilde{\bp},\tilde{\bla})^*$. Now all $\tilde{p}_i$ is odd, it follows that $B(\tilde{\bp},\tilde{\bla})$ and $\tilde{B}(\bp,\bla)$ is invertible.

    \bexample
    Let $N=3$, $\bp=(1,1,4)$ and $\bla=(0,\infty,1)$, then the Grothendieck group $K_0$ of the coherent category $\mathrm{Coh}(\X_{\bp,\bla})$ has a basis
         $$\{\hat{\sO}, \hat{S_x},\hat{S}_{1,2}, \hat{S}_{1,3},\hat{S}_{1,4}\}.$$
    Therefore $B(\bp,\bla)$ looks like
        $$\begin{bmatrix}
        0 &-2 &0 &0 &-1\\
        2 &0   &0  &0 &0\\
        0 &0   &0  &1 &0  \\
        0 &0  &-1 &0 &1 \\
        1 &0    &0    &-1 &0 \\
         \end{bmatrix}.$$
        Set $\tilde{\bp}=(1,1,5)$, then $B(\tilde{\bp},\bla)$ is as follows:
        $$\begin{bmatrix}
          0 &-2 &0 &0 &0  &-1\\
          2 &0   &0  &0 &0 &0\\
          0 &0   &0  &1 &0  &0 \\
          0 &0  &-1 &0 &1 &0 \\
          0 &0   &0  &-1 &0 &1\\
          1 &0 &0 &0 &-1 &0
           \end{bmatrix}.$$
        Move the third column to the last and then the third row to the last, we get $\tilde{B}(\bp,\bla)$:
        $$\begin{bmatrix}
          0 &-2  &0 &0  &-1 &0\\
          2 &0    &0 &0 &0  &0\\
          0 &0    &0 &1 &0  &-1 \\
          0 &0     &-1 &0 &1 &0\\
          1 &0   &0 &-1 &0  &0\\
          0 &0     &1 &0  &0  &0\\
           \end{bmatrix}.$$
      It can be easily checked that $\tilde{B}(\bp,\bla)$ is invertible and $B(\bp,\bla)$ is the upper-left submatrix of $\tilde{B}(\bp,\bla)$.
    \eexample

    Fix $\tilde{B}(\bp,\bla)$ constructed as above. Since $\tilde{B}(\bp,\bla)$ is skew-symmetric and invertible, there exists $d\in \N$ and an $m\times m$ skew-symmetric matrix $\varLambda$ of integers such that
                    $$-\varLambda \tilde{B}(\bp,\bla)=d I_m.$$

    \bremark
    In general, $d$ may not be 1. In the level of categorifications, it can be realized by working on everything over the field $\F_{q^d}$ rather than $\F_q$.  Let $\tilde{B}(\bp,\bla)$ be constructed as above. Say  $-\varLambda' \tilde{B}(\bp,\bla)=I_m$ for some skew-symmetric matrix $\varLambda'$ such that $d\varLambda'$ is a matrix of integers. Then $|\F_{q^d}|^{\frac{1}{2}\varLambda'(x,y)}=q^{\frac{d}{2}\varLambda'(x,y)}$ will be a polynomial of $v^{\pm1}$, where $v=q^{\frac{1}{2}}$.
    \eremark
  
    In the sequel, we will write $\tilde{E}'$ for $\tilde{E}(\tilde{\bp},\tilde{\bla})^{*t}$, and $\tilde{E}$ for $\tilde{E}(\tilde{\bp},\tilde{\bla})^*$. 
    
  \bproposition\label{prop1}
  We identify $\varLambda$ with bilinear form and have the following 
   We have the following identities for $\bm,\bn\in K_0(\tilde{\sA})$:
   \begin{enumerate}
   \item[(i)]$\varLambda(\tilde{B}\bm,\tilde{E}\bn)= \langle \bm,\bn \rangle $.

   \item[(ii)]$\varLambda(\tilde{B}\bm,\tilde{E}'\bn)= \langle \bn,\bm \rangle $.
   
   \item[(iii)]$\varLambda(\tilde{B}\bm,\tilde{B}\bn)= \langle \bn,\bm \rangle - \langle \bm,\bn \rangle $.

   \item[(iv)]$\varLambda(\tilde{E}\bm,\tilde{E}\bn)=\varLambda(\tilde{E}'\bm,\tilde{E}'\bn)$.
   \end{enumerate}
   \eproposition

   \bproof
    These results are obtained by direct computation.
   \eproof

   Denote $\bm^*:=\tilde{E}'\bm$ and $^*\bn:=\tilde{E}\bn$.
   
   \blemma
     $\varLambda(-b^*-^*a,-d^*-^*c)=\varLambda((a+b)^*,(c+d)^*)+ \langle b,c \rangle - \langle d,a \rangle .$
   \bproof
    we have
    \begin{align*}
      &\varLambda(-b^*-^*a,-d^*-^*c)\\
      &=\varLambda(b^*+a^*,d^*+c^*)-\varLambda(b^*,c^*)-\varLambda(a^*,d^*)+\varLambda(b^*,^*c)+\varLambda(^*a,d^*)\\
      &=\varLambda(b^*+a^*,d^*+c^*)+ \langle b,c \rangle - \langle d,a \rangle .
    \end{align*}
    The first equality is induced by $(iv)$ and the second is by $(i)$ and $(ii)$ in Proposition \ref{prop1}.
   \eproof
   \elemma

   \subsection{Quantum torus and integration maps}
    Let $B(\bp,\bla)$ be the skew-symmetric matrix associated to $\sA$. Then there exist $(\tilde{\bp},\tilde{\bla})$ and a skew-symmetric matrix $\varLambda$ such that $-\varLambda B(\tilde{\bp},\tilde{\bla})^*=I_m$.  Notice that there exists some positive integer $d$ such that $d \varLambda$ is a matrix of integers. We take $d$ as the minimal one. Denoted by $\tilde{\sA}_k$  the category $\mathrm{Coh}(\X_{\tilde{\bp},\tilde{\bla}})_k$ over $k=\F_{q^d}$.
   
    let $\nu$ be a formal invariable. $\sT_m$ is defined to be the $\Z[\nu^{\pm1}]$-algebra with a basis $\{X^{\alpha}|\alpha\in \Z^m\}$ (namely, $\sT_m=\Z[\nu^{\pm1}][x_1^{\pm1},x_2^{\pm1},\cdots,x_m^{\pm1}]$, where $x_i$ are formal variables) and multiplication given by 
    $$X^\alpha X^{\beta}=X^{\alpha+\beta}.$$
     The quantum torus $\sT_\varLambda$ is a  $\Z[\nu,\nu^{-1}]$-algebra with the same vector space as $\sT_m$ but with a twisted multiplication:
    $$X^\alpha * X^\beta=\nu^{d\varLambda(\alpha^*,\beta^*)}X^{\alpha+\beta}.$$

   Set $v=q^{\frac{1}{2}}$.  Denoted by $\sT_{\varLambda,v}$ (resp. $\sT_{\varLambda}$) the specialization of $\sT_m$ (resp. $\sT_{\varLambda}$) at $\nu=v$. Let $\hat{\sT}_{\varLambda}$ (resp. $\hat{\sT}_{\varLambda,v}$) be the completion 
    $$\Z[\nu^{\pm1}][X^{-\bff_1},X^{\bff_2},X^{-\bff_3},\cdots,X^{-\bff_m}][[X^{\bff_1},X^{-\bff_2},X^{\bff_3},\cdots,X^{\bff_m}]],$$
   of $\sT_{\varLambda}$ (resp. $\sT_{\varLambda,v}$), where $f_i=\tilde{B}(\bp,\bla)\bee_i$.
    
   \begin{proposition}
   The integration map $\int: H^{\vee}(\tilde{\sA}_{k})\lrw \sT_{m,v}$, $[\sF]\mapsto X^{\underline{\dim}\sF}$ is an algebraic homomorphism.
    \begin{proof}
      We have
    \begin{align*}
       \int([\sF][\sG])  &=\sum_{[\sL]}q^{-d[\sF,\sG]^1} |\operatorname{Ext}^1_{\tilde{\sA}}(\sF,\sG)_{\sL}| X^{\underline{\dim}\sL}\\
                                &=q^{-d[\sF,\sG]^1} |\operatorname{Ext}^1_{\tilde{\sA}}(\sF,\sG)| X^{\underline{\dim} \sF+\underline{}{\dim}\sG}\\
                                &=X^{\underline{\dim}\sF}X^{\underline{\dim}\sG}\\
                                &=\int([\sF])\int([\sG]).
    \end{align*}
    \end{proof}
   \end{proposition}
  
   \subsection{\texorpdfstring{$\varLambda$}{Lg}-twisted versions}
   Provided the skew-symmetric form $\varLambda$,  we twist the multiplication on $H^{\vee}(\tilde{\sA})$ as follows:
        $$[\sM]*[\sN]=v^{d\varLambda(m^*,n^*)}[\sM][\sN].$$
    where $\bm$ (resp. $\bn$) is the dimension vector of $[\sM]$ (resp.$[\sN]$) in $K_0(\tilde{\sA}_{k})$. The $\varLambda$-twisted Hall algebra is denoted by $H_{\varLambda}(\tilde{\sA}_{k})$.

   We also twist the multiplication on $H^{\vee}(\tilde{\sA}_{k})\hat{\otimes} H^{\vee}(\tilde{\sA}_{k})$ again such that the coproduct $\Delta$ is still an algebra homomorphism.  Let $(H^{\vee}(\tilde{\sA}_{k})\hat{\otimes} H^{\vee}(\tilde{\sA}_{k}),*)$ be the tensor algebra with twisted multiplication $*$ given as
            $$([\sM_1]\otimes [\sM_2])*([\sN_1]\otimes[\sN_2]):=v^{d\varLambda((\bm_1+\bm_2)^*,(\bn_1+\bn_2)^*)}([\sM_1]\otimes [\sM_2])([\sN_1]\otimes[\sN_2]).$$
   Hence it can be easily checked that $\Delta: H_{\varLambda}(\tilde{\sA}_{k})\to (H^\vee(\sA_{k})\hat{\otimes} H^{\vee}(\sA_{k}),*)$ is also an algebra homomorphism.

   Recall that we have defined an integration map $\int: H^{\vee}(\tilde{\sA}_{k})\to \sT_{m,v}$, which induces a map 
          $$\int\otimes \int:H^{\vee}(\tilde{\sA}_{k})\hat{\otimes} H^{\vee}(\tilde{\sA}_{k})\to \sT_{m,v}\hat{\otimes} \sT_{m,v},\ [\sM]\otimes [\sN]\mapsto X^{\bm}\otimes X^{\bn}.$$
   Note if the multiplications on $H^{\vee}(\tilde{\sA}_{k})\hat{\otimes} H^{\vee}(\tilde{\sA}_{k})$ and $\sT_{m,v}\hat{\otimes} \sT_{m,v}$ both are untwisted, that is $(x_1\otimes y_1)(x_2\otimes y_2)=x_1x_2\otimes y_1y_2$, then $\int\otimes \int$ is a homomorphism of algebras. Since we have twisted the multiplication on $H^{\vee}(\tilde{\sA}_{k})\hat{\otimes} H^{\vee}(\tilde{\sA}_{k})$, we also twisted the multiplication on $\sT_{m,v}\hat{\otimes} \sT_{m,v}$ by 
              $$(X^{\alpha_1}\otimes X^{\beta_1})*(X^{\alpha_2}\otimes X^{\beta_2}):=q^{\frac{d}{2}\varLambda((\alpha_1+\beta_1)^*,(\alpha_2+\beta_2)^*)+d(\beta_1,\alpha_2)+ d\langle \alpha_1,\beta_2 \rangle }X^{\alpha_1+\alpha_2}\otimes X^{\beta_1+\beta_2}.$$
   Then $\int\otimes \int: (H^{\vee}(\tilde{\sA}_{k})\hat{\otimes} H^{\vee}(\tilde{\sA}_{k}),*)\lrw (\sT_{m,v}\hat{\otimes} \sT_{m,v},*)$ is a homomorphism of algebras.

   Finally, following \cite[Proposition 7.11]{Fu2020}, we define an algebra homomorphism $\mu: (\sT_{m,v}\hat{\otimes} \sT_{m,v},*)\lrw \hat{\sT}_{\varLambda,v}$ by
    $$\mu(X^\alpha\otimes X^\beta)=v^{-d(\alpha,\beta)- d\langle \alpha,\beta \rangle }X^{-^{*}\alpha-\beta^*}.$$
   
   Therefore we get an algebra homomorphism $X_?: H_{\varLambda}(\tilde{\sA}_{k})\lrw \hat{\sT}_{\varLambda,v}$ given by the composition $\mu\circ(\int\otimes\int)\circ\Delta$, called the character map. Namely, we have the following commutative diagram:
    $$\begin{tikzcd}
      H_{\varLambda}(\tilde{\sA}_{k})\arrow[r,"X_?"]\arrow[d,"Delta"]    &\hat{\sT}_{\varLambda,v}\\
     (H^{\vee}(\tilde{\sA}_{k})\hat{\otimes} H^{\vee}(\tilde{\sA}_{k}),*)\arrow[r,"\int\otimes \int"]     &(\sT_{m,v}\otimes\sT_{m,v},*)\arrow[u,"\mu"].\\
    \end{tikzcd}$$

    Therefore, for $\sM\in \tilde{\sA}_{k}$, the quantum cluster character of $\sM$ is 
          \begin{align*}
              X_{\sM}&=\sum_{\sX,\sY} q^{-\frac{d}{2} \langle \sY,\sX \rangle }g_{\sX,\sY}^{\sM} X^{-y^*-^*x}\\
                 &=\sum_{\boldsymbol{e}\leq \underline{\dim}\sM} q^{-\frac{d}{2} \langle \bm-\boldsymbol{e},\boldsymbol{e} \rangle }|\mathrm{Gr}_{\boldsymbol{e}}(\sM)| X^{-(\bm-\boldsymbol{e})^*-^*\boldsymbol{e}}.
          \end{align*}
    where $y=\underline{\dim}\sF$,  $x=\underline{\dim}\sG$ and $\bm=\underline{\dim}\sM$, $\mathrm{Gr}_{\boldsymbol{e}}(\sM)$ is the Grassmannian variety of subobjects of $\sM$ with dimension vector $\boldsymbol{e}$ and $|\mathrm{Gr}_{\bee}(\sM)|$ is its cardinality.

  \bexample\label{ex2.5} Let $B$ be the skew-symmetric matrix associated to $\P^1$. Then 
   $$B=\begin{bmatrix}0  &-2\\ 2  &0\end{bmatrix}=\begin{bmatrix}1  &-1\\ 1  &0\end{bmatrix}-\begin{bmatrix}1  &1\\ -1  &0\end{bmatrix},$$
   $\varLambda=\begin{bmatrix}0  &-\frac{1}{2}\\ \frac{1}{2}  &0\end{bmatrix}$ and $d=2$. Hence $\tilde{\sA}$ is the category $\mathrm{Coh}(\P^1)$ over $\F_{q^2}$. 

  (1) $ X_{\sO(l)}=X^{-(\begin{smallmatrix}1+l\\ -1\end{smallmatrix})}+\sum_{r\geq -l} v^{-2(l+r)}[l+r+1]_{q'} X^{-(\begin{smallmatrix}l+2r+1\\ 1\end{smallmatrix})}$, here $[n]_{q}$ means $\frac{q^{n}-1}{q-1}$, and $q'=q^2$.

  (2) Let $S_x$ be a simple torsion sheaf supported on $x\in \P^1$ with degree $d$, then 
      $$X_{S_x}=X^{-(\begin{smallmatrix}d\\ 0\end{smallmatrix})}+X^{-(\begin{smallmatrix}-d\\ 0\end{smallmatrix})}.$$
  \eexample

  \subsection{Definition of the quantum cluster algebra of \texorpdfstring{$\X_{\bp,\bla}$}{Lg}}\label{sec3.4}
   In this subsection, we want to define the quantum cluster algebra of the weighted projective line $\X_{\bp,\bla}$. 
   
   Recall the definition of quantum cluster algebras  introduced by \cite{Berenstein2005}. Let $n\leq m$, ${\sT}_{\varLambda}={\sT}(\Z^m,\varLambda)$be the quantum torus.  Let $(\varLambda,\tilde{B},X)$ be an initial seed (see \cite[Definition 2.1.5]{Qin2012}), and $\T_m$ be an $m$-regular tree with root $t_0$. By \cite[Corollary 2.1.10]{Qin2012},  given another seed $(\varLambda',\tilde{B}',X')$, we say $X'$ is mutated from $X$ at $i$ $(1\leq i\leq n)$ if 
   \begin{itemize}
    \item[$\cdot$] $X(e_j)=X'(e_j)\qquad$ if $j\neq i$,\\
    \item[$\cdot$]  \begin{equation}\label{eq1}
      \begin{aligned}
      &X(e_i)X'(e_i)\\
      &=v^{\varLambda(e_i,\sum\limits_{1\leq l\leq m}[b_{li}]_+ e_l)}X({\sum\limits_{1\leq l\leq m}[b_{li}]_+ e_l})+v^{\varLambda(e_i,\sum\limits_{1\leq l\leq m}[-b_{li}]_+ e_l)}X({\sum\limits_{1\leq l\leq m}[-b_{li}]_+ e_l}).
    \end{aligned}
  \end{equation}
   \end{itemize}
   
   Write $\begin{tikzcd}t\arrow[r,no head] &t'\end{tikzcd}$ if $t$ and $t'$ of $\T_m$ are linked by an edge labeled $i$. Then one can associate  iteratively each seed mutated from $X$ with each vertex $t$ of $\T_m$. Namely,  Set the initial seed to be $(\varLambda(t_0),\tilde{B}(t_0),X(t_0))$. If $\begin{tikzcd}t\arrow[r,no head] &t'\end{tikzcd}$ and $1\leq i\leq n$, then label the seed mutated from $(\varLambda(t),\tilde{B}(t),X(t))$ at $i$ by $(\varLambda(t'),\tilde{B}(t'),X(t'))$;  If $\begin{tikzcd}t\arrow[r,no head] &t'\end{tikzcd}$ and $n<i\leq m$, then set $(\varLambda(t'),\tilde{B}(t'),X(t'))=(\varLambda(t),\tilde{B}(t),X(t))$. The quantum cluster algebra of $(\varLambda,\tilde{B},X)$ is defined to be a $\Z[\nu^{\pm 1}]$-subalgebra of the quantum torus $\sT_{\varLambda}$ generated by quantum cluster variables $X_i(t)$ for all the vertices $t\in \T_m$, $1\leq i\leq n$ and elements $X_j(t_0)^{\pm1}$ for all $n<j\leq m$ in \cite{Qin2012}. There Qin used the refined CC-map to categorify the quantum cluster algebra, and showed that the CC-map $X_{T_i(t)}$ of indecomposable coefficient-free rigid objects $T_i(t)$ in certain cluster category are bijectively corresponding to the quantum cluster variables $X_i(t)$ for $t\in \T_m$. In a similar way, we will give the definition of quantum cluster algebra of a weighted projective line by setting generators indexed by certain indecomposable rigid objects in a cluster category.
   
   For a field k, let $\tilde{\sA}_k$ be the category of coherent sheaves on $\X_{\tilde{\bp},\tilde{\bla}}$ over k. The cluster category $\sC_k:=\sC(\tilde{\sA}_k)$ is defined to be the orbit category $D^b(\tilde{\sA}_k)/\tau\circ[-1]$, where $\tau$ is the Auslander-Reiten translation. Following \cite[Theorem 6.8]{Buan2006a}, any almost cluster-tilting object $\bar{T}^k$ has exactly two complements $T_i^k$ and $T_i^{*k}$. Such $(T_i^k,T_i^{*k})$ is called an exchange pair. Moreover, $T_i^k$ and $T_i^{*k}$ are linked by exchange triangles:
      $$T_i^k\stackrel{u}\lrw E^k \stackrel{v}\lrw T_i^{*k} \lrw T_i^k[1],\qquad\text{and}\qquad 
      T_i^{*k}\stackrel{u'}\lrw E^{'k} \stackrel{v'}\lrw T_i \lrw T^{*k}_i[1],$$
   where $u$ and $u'$ are minimal left $\mathrm{add}\bar{T}^k$-approximations and $v$ and $v'$ are minimal right $\mathrm{add}\bar{T}^k$-approximations. Write $E^k=\bigoplus_{j \neq i} T_j^{k,\oplus a_{ij}}$ and $E^{'k}=\bigoplus_{j\neq i} T_j^{k,\oplus b_{ij}}$. Let $A_k:=\mathrm{End}_{\sC_k}(T^k)$ and $Q_{T^k}$ be the Gabriel quiver of $A_k$. Since we have an equivalence $\mathrm{add}\bar{T}^k \stackrel{\sim}\lrw \mathrm{proj}A_k$ of additive categories, $a_{ij}(k)$ is the number of arrows from $j$ to $i$ in $Q_{T^k}$ and $b_{ij}(k)$ is the number of arrows from $i$ to $j$ in $Q_{T^k}$.  Hence, we can also construct a $m$-regular tree $\T_m(k)$ as above, where $m$ is the rank of $K_0(\tilde{\sA}_k)$. Furthermore. we need to record the number of arrows $(a_{ij}(k),b_{ij}(k))$ of the Gabriel quiver $Q_{T^k}$ of the endomorphism algebra of each cluster-tilting object $T^k$ in $\sC_k$.  Let $T^k:=\bigoplus_{0\leq \vec{l}\leq \vec{c}} \sO(\vec{l})^k$ be an initial cluster-tilting object in $\sC_k$, which is associated to the root $t_0$ of the tree $\T_m(k)$. i.e. $T(t_0)=T^k$. If $T^{'k}$ is mutated from $T^k$ at the $i$-th direct summand $T_i^k$, we set $T(t)=T^{'k}$ where $\begin{tikzcd}[column sep=huge]t_0\arrow[r, no head,"{(a_{ij}(k),b_{ij}(k))}"] &t\end{tikzcd}$. Here $a_{ij}(k)$ (resp. $b_{ij}(k)$) is the number of arrows from $j$ to $i$ (resp. $i$ to $j$) of the quiver $Q_{T(t_0)}$.

   \bdefinition\label{def3.6}
   The regular m-tree $\T_m(T^k)$ constructed as above is called the valued regular $m$-tree over $k$ with the initial cluster-tilting object $T^k$ associated to $\sC(\tilde{\sA}_k)$.
   \edefinition
  
   According to Theorem \ref{thm7.8}, we know that each valued regular m-tree $\T_m(T^{\F_{q^r}})$ over $\F_{q^r}$ is the same as $\T_m(T^{\bar{\F}_q})$ for a fixed prime $q$ and  $r\geq 1$. To show that $\T_m(T^{\bar{\F}_q})$ and $\T_m(T^{\bar{\F}_p})$ are identical for distinct primes $p$ and $q$, we need the following 
  \btheorem[{\cite[Theorem 5.2]{Buan2011}}]
   Let $\sC_K$ be a 2-$CY$ triangulated category with a cluster-tilting object $T$ over an algebraically closed field $K$. If the endomorphism algebra $\mathrm{End}_{\sC_K}(T)$ is isomorphic to the Jacobian algebra $J(Q,W)$ for some quiver with potential $(Q,W)$, and if no 2-cycles start in the vertex $i$ of $Q$, then we have an isomorphism
      $$\mathrm{End}_{\sC_K}(\mu_i(T))=J(\mu_i(Q,W))  .$$
   \etheorem
   Here $\mu_i(Q,W)$ is the mutation of quiver with potentials, see \cite[Section 1.2]{Buan2011}. Since $\tilde{\sA}_K$ is derived equivalent to a canonical algebra which is of global dimension $\leq 2$, by \cite[Theorem 6.12]{Keller2011} $\mathrm{End}_{\sC_K}(T)=J(Q,W)$ for some quiver with potential $(Q,W)$, where $T:=\bigoplus_{0\leq \vec{l}\leq \vec{c}} \sO(\vec{l})$. Moreover, the quiver with potential $(Q,W)$ is non-degenerate by \cite[Lemma 3.2]{Fu2021}. Combining with the theorem above and notice that $Q_{T^{\bar{\F}_q}}$ is the same as $Q_{T^{\bar{\F}_p}}$, we can conclude that the quiver $Q_{T^{'K}}$ of each cluster-tilting object $T^{'K}$ mutated from $T^K$ is independent of the choice of algebraically closed fields. Thus $\T_m(T^{\bar{\F}_q})=\T_m(T^{\bar{\F}_q})$ for any primes $p$ and $q$. On the other hand, because the cluster-tilting graph of $\sC(\tilde{\sA}_k)$ is connected by Corollary \ref{cor7.9}, the valued regular $m$-tree $\T_m(T^k)$ is also independent of the choice of the initial cluster-tilting objects $T^k$. Namely, if $T^{'k}$ is another cluster-tilting object, then $T^{'k}=T(t)$ for some $t\in \T_m(T^k)$ and the valued regular $m$-tree $\T_m(T^{'k})$ with initial object $T^{'k}$ is obtained from $\T_m(T^k)$ by taking $t$ to be the new root. In summary, we obtain the following
   \blemma
   The valued regular m-tree $\T_m(T^k)$ associated to $\sC(\tilde{\sA}_k)$ with initial cluster-tilting object $T^k$ is independent of finite fields and the choice of initial cluster-tilting objects.
   \elemma
   Set the initial cluster-tilting object to be
     $$S^k:=\tilde{\sO}^k\oplus \tilde{\sO}(\vec{\tilde{c}})^k\oplus \bigoplus_{2\leq j\leq \tilde{p}_i} \tilde{S}_{i,j}^k.$$ 
  In the sequel, we will abbreviate $\T_m$ for $\T_m(S^k)=\T_m(T^k)$.
   
   One of exchange triangles linking $T_i^k$ and $T_i^{*k}$ is induced by a short exact sequence in $\tilde{\sA}_k$, say $\mathrm{Ext}^1_{\tilde{\sA}_k}(T_i^{*k}, T_i^k)\cong k$ and the other is given by
           $$T_i^{*k}\lrw E^{'k}\lrw T_i^k\stackrel{f_k}\lrw T_i^{*k}[1],$$
    where $E^{'k}\cong \Ker f_k\oplus \tau^{-1}\mathrm{Coker} f_k$. Hence in $K_0(\tilde{\sA}_k)$ we have $[T_i^{*k}]=[E^k]-[T_i^k]=[E^{'k}]-[T_i^k]+[\mathrm{Im} f_k]+[\tau^{-1}\mathrm{Im} f_k]$, which implies that $[E^k]-[E^{'k}]=[\mathrm{Im} f_k]+[\tau^{-1}\mathrm{Im} f_k]>0$. Similarly, if $\mathrm{Ext}^1_{\tilde{\sA}_k}(T_i^k,T_i^{*k})\cong k$, then we have $[E^{'k}]-[E^k]=[\mathrm{Im} g_k]+[\tau^{-1}\mathrm{Im} g_k]>0$, where $g_k: T_i^k\to T_i^{*k}[1]$ is nonzero. Note that if dimension vectors of $[E^k]$, $[E^{'k}]$ and $[T_i^k]\in K_0(\tilde{\sA}_k)$ are independent of the choice of fields, then for two distinct fields $k_1$ and $k_2$, we have $\mathrm{Ext}^1_{\tilde{\sA}_{k_1}}(T_i^{*k_1}, T_i^{k_1})\cong k_1$ implying that $[E^{'k_2}]-[E^{k_2}]=[E^{'k_1}]-[E^{k_1}]<0$ and then $\mathrm{Ext}^1_{\tilde{\sA}_{k_2}}(T_i^{k_2},T_i^{*k_2})$ must be $0$. As a consequence, 
    $$(*)\qquad \dim_{k_1}\mathrm{Ext}^1_{\tilde{\sA}_{k_1}}(T_i^{*k_1},T_i^{k_1})= 1 \text{ if and only if } \dim_{k_2}\mathrm{Ext}^1_{\tilde{\sA}_{k_2}}(T_i^{*k_2},T_i^{k_2})= 1.$$ 
   Hence, by induction from the root $t_0$ of the tree $\T_m(k)$, it can be showed that the dimension vector of $[T_i(t)]$ for $t\in \T_m$ is independent of the choice of fields. So we will write $d_i(t)$ for the dimension vector $[T_i(t)^k]$ for $t\in \T_m$, $1\leq i\leq m$. 
  
    Let $B(\bp,\bla)$ be the skew-symmetric matrix associated to $\sA$, and $(\varLambda, \tilde{B}(\bp,\bla))$ the compatible pair given as in Section \ref{sec3.1}, where $\tilde{B}:=\tilde{B}(\bp,\bla)$ is similar to $B(\tilde{\bp},\tilde{\bla})$.   If $B(\bp,\bla)$ is a proper submatrix of $\tilde{B}(\bp,\bla)$. i.e. $m>n$, then it does not need to do mutations at every direct summand of the initial cluster-tilting object $S$. If $\tilde{\bp}=(p_1+1,p_2,\cdots,p_N)$, then $\phi_*^{-1}(S')=S/\tilde{S}_{1,2}$ where $S'=\sO\oplus \sO(\vec{c})\oplus \bigoplus_{2\leq j\leq p_i} S_{i,j}$ is a cluster-tilting object in $\sC(\sA)$. Therefore, it does not need to do mutations at $\tilde{S}_{i,2}$ for $p_i$ even. Write $T(t_0)=\bigoplus_{i=1}^m T_i(t_0)=S$ , we order the direct summands of $S$ as the basis $\boldsymbol{b}(\tilde{\bp},\tilde{\bla})^*$ for $K_0(\tilde{\sA})$ defined in Section \ref{sec3.1}. Then, we may not do mutations at $T_j(t)^k$ of the cluster-tilting object $T(t)^k$ for $n<j\leq m$, $t\in \T_m$, it follows that the subgraph of $\T_m$ consisting of vertices where we actually do mutations with respect to $(\bp,\bla)$ is a regular $n$-tree, denoted by $\T_n(\bp,\bla)$.

    In the sequel, we let $T_i(t)$ and $T(t)$ be symbols associated to $t\in \T_m$, $1\leq i\leq m$. Note that $T_i(t)^k$ (resp. $T(t)$) is a indecomposable rigid (resp. cluster-tilting) object in $\sC(\tilde{\sA}_k)$ labeled by $t\in \T_m$ for some $i$. By $(*)$, $\mathrm{Ext}^1(T_i(t),T_i(t'))=0$ means that $\mathrm{Ext}^1_{\tilde{\sA}_k}(T_i(t)^k,T_i(t')^k)\neq 0$.

    Now we are in the position to give the definition of quantum cluster algebra of $\X_{\bp,\bla}$.  For convenience of notations, we assume that $-\varLambda \tilde{B}(\bp,\bla)=I_m$ with $\varLambda\in \mathrm{Mat}(m,\Z)$. If $-\varLambda' \tilde{B}(\bp,\bla)=I_m$ such that $d\varLambda'$ is a matrix of integers for some $d$, then we only need to set the following $\nu$ to be $\nu^d$.

    Recall 
      $$S^k=\tilde{\sO}^k\oplus \tilde{\sO}(\vec{\tilde{c}})^k\oplus \bigoplus_{2\leq j\leq \tilde{p}_i} \tilde{S}_{i,j}^k=:\bigoplus_{i=1}^m T_i(t_0)^k$$
     is the initial cluster-titling object in $\sC(\tilde{\sA}_k)$. Every subject of the line bundle $\sO(\vec{l})^k$ is of the form $\sO(\vec{r})^k$ such that $\vec{r}=\sum_{i=1}^N r_i\vec{x_i}+r_0\vec{c}\leq \vec{l}=\sum_{i=1}^N l_i\vec{x_i}+l_0\vec{c}$. Hence, for each $\bee\leq \dim[\sO(\vec{l})]$, there exists a unique isoclasses $[\sO(\vec{r_e})^k]$ with dimension vector $\bee$ such that $\sO(\vec{r_e})^k$ is a subject of $\sO(\vec{l})^k$.  It is easy to see that there exists a $\Z$-polynomial $P(z)$ (independent of $k$) such that $|\mathrm{Gr}_{\bee}(\sO(\vec{l})^k)|=P_i(|k|)$. Denote $|\mathrm{Gr}_{\bee}(T_i(t_0))|_{\nu^2}:=P_i(\nu^2)$ for $i=1,2$. Set
       $$  X_i(t_0)=\sum_{\boldsymbol{e}\leq d_i(t_0)} \nu^{\langle d_i(t_0)-\boldsymbol{e},\boldsymbol{e} \rangle }|\mathrm{Gr}_{\boldsymbol{e}}(T_i(t_0))|_{\nu^2} X^{-(d_i(t_0)-\boldsymbol{e})^*-^*\boldsymbol{e}}$$
    for $i=1,2$. For $3\leq i\leq m$, $T_i(t_0)^k=\tilde{S}_{l,j}^k$ is a simple torsion sheaf for some $l,j$, then we set
       $$X_i(t_0)=X^{-\bee_i^*}(X^{\tilde{B}\bee_i}+1),$$
    where $\{\bee_i\ | 1\leq i \leq m\}$ is the canonical basis for $\Z^m$. Observe that $X_i(t_0)=x_l+1$  for $n<i\leq m$, since $\tilde{B}\bee_i=-\bee_l$ where $3\leq l\leq m$ such that $\bee_l=\underline{\dim}[\tilde{S}_{j,3}]$. Here $T_i(t_0)=\tilde{S}_{j,2}$ for some $j\in \{1\leq i\leq N|\ p_i \text{ is even}\}$. Hence $X_i(t_0)$ is invertible in $\hat{\sT}_{\varLambda}$ for $i>n$. 
   
    We define another partial order in $\Z^m$ associated to $\bff\in \Z^m$. Notice that $\tilde{B}(\boldsymbol{e})=\boldsymbol{e}^*-^*\boldsymbol{e}$ and $\tilde{B}$ is invertible, then $\bm=\tilde{B}\bee$ is  uniquely determined by $\bee$. We say 
         $$ -^*\bff+\tilde{B}\bee\leq -^*\bff+\tilde{B}\bee' \text{    if and only if    } \bee\leq \bee'.$$
    Hence the maximal degree of $X_i(t_0)$ is $-^*d_i(t_0)$.

    \blemma\label{lem3.9}
    $\{X_i(t_0)|1\leq i\leq m\}$ is algebraically independent in $\hat{\sT}_{\varLambda}$.
    \bproof
    $T(t_0)^k=\tilde{\sO}^k\oplus \tilde{\sO}(\vec{\tilde{c}})^k\oplus \bigoplus_{2\leq j\leq \tilde{p}_i} \tilde{S}_{i,j}^k$ is a cluster-tilting object, the set $\{d_i(t_0)| 1\leq i\leq m\}$ forms a basis for $K_0(\tilde{\sA})$. Hence the maximal degrees $-^*d_i(t_0)$ of $X_i(t_0)$, $1\leq i\leq m$ forms a basis of $\Z^m$ by noting that  $E(\tilde{\bp},\tilde{\bla})$ is invertible. Since 
     $$\{X^{\bee_i}=x_i\ |\ 1\leq i\leq m\}$$
    are algebraically independent in $\hat{\sT}_{\varLambda}$, it follows that $\{X^{-^*d_i(t_0)}|1\leq i\leq m\}$ are algebraically independent. As a consequence,  $\{X_i(t_0)\ |\ 1\leq i\leq m\}$ is algebraically independent in $\hat{\sT}_{\varLambda}$.
    \eproof
    \elemma

   \bdefinition\label{def3.10}
   The quantum cluster algebra $\sA(\varLambda,\tilde{B}(\bp,\bla))$ of the weighted projective line $\X_{\bp,\bla}$ is the $\Z[\nu^{\pm 1}]$-subalgebra of $\hat{\sT}_{\varLambda}$, generated by $X_j(t)$ for $t\in \T_n(\bp,\bla)$, $1\leq j\leq n$ and $X_l(t_0)^{\pm 1}$ for $n<l\leq m$, subject to 
   \begin{itemize}
    \item[(1)] if $j\neq i$
      $$X_j(t)X_i(t)=\nu^{2\varLambda(d_j(t)^*,d_i(t)^*)}X_i(t)X_j(t),$$
    \item [(2)] if  $\begin{tikzcd}[column sep=huge]t\arrow[r, no head,"{(a_{ij}(t),b_{ij}(t))}"] &t'\end{tikzcd}$ and $\mathrm{Ext}^1(T_i(t),T_i(t'))=0$,
    $$X_i(t')X_i(t)=\nu^{\varLambda(d_i(t')^*,d_i(t)^*)} \nu^s \prod_{j\neq i}X_j(t)^{a_{ij}(t)}  
    + \nu^{\varLambda(d_i(t')^*,d_i(t)^*)-1} \nu^{s'} \prod_{j\neq i}X_j(t)^{b_{ij}(t)},$$
    \item [(3)] if  $\begin{tikzcd}[column sep=huge]t\arrow[r, no head,"{(a_{ij}(t),b_{ij}(t))}"] &t'\end{tikzcd}$ and $\mathrm{Ext}^1(T_i(t'),T_i(t))=0$,
    $$X_i(t')X_i(t)=\nu^{\varLambda(d_i(t')^*,d_i(t)^*)+1} \nu^s \prod_{j\neq i}X_j(t)^{a_{ij}(t)}  
    + \nu^{\varLambda(d_i(t')^*,d_i(t)^*)} \nu^{s'} \prod_{j\neq i}X_j(t)^{b_{ij}(t)},$$
   \end{itemize}

   where $s=-\sum\limits_{l=1}^m \varLambda(a_{il}d_l(t)^*,\sum\limits_{r=l+1}^n a_{ir}d_r(t)^*)$, $s'=-\sum\limits_{l=1}^m \varLambda(b_{il}d_l(t)^*,\sum\limits_{r=l+1}^n b_{ir}d_r(t)^*)$.
   \edefinition

   \bremark\label{rem2}\ \ 

   (1) Although the mutation relations in Definition \ref{def3.10} are similar to Relations (\ref{eq1}) of the usual quantum cluster algebra, the exchange matrix $\begin{tikzcd}[column sep=huge]t_0\arrow[r, no head,"{(a_{ij}(t_0),b_{ij}(t_0))}"] &t\end{tikzcd}$ has nothing to do with $\tilde{B}(\bp,\bla)$. Indeed, $(a_{ij}(t_0),b_{ij}(t_0))$ can be read from the Gabriel quiver $Q_{T_0}$ of $\mathrm{End}_{\sC}(T_0)$, which is (see \cite[Theorem 6.12]{Keller2011})
    $$\begin{tikzcd}
       &\circ\arrow[ldd,"\alpha_{1,1}"]  &\circ\arrow[l,"\alpha_{1,2}"] &\cdots\arrow[l,"\alpha_{1,p_1-2}"] &\circ\arrow[l,"\alpha_{1,p_1-1}"] & \\
       & &\vdots & & &\\
      \star\arrow[rrrrr,"\rho_3",bend left]\arrow[rrrrr,"\rho_N",bend right]\arrow[rrrrr,"\rho_4",bend left=20] \ &\circ\arrow[l,"\alpha_{i,1}"] &\circ\arrow[l,"\alpha_{i,2}"] &\cdots\arrow[l,"\alpha_{i,p_i-2}"] &\circ\arrow[l,"\alpha_{i,p_i-1}"] &\ast\arrow[luu,"\alpha_{1,p_1}"]\arrow[l,"\alpha_{i,p_i}"]\arrow[ldd,"\alpha_{N,p_N}"']   \\
       & &\vdots  & & &\\
       &\circ\arrow[luu,"\alpha_{N,1}"']  &\circ\arrow[l,"\alpha_{N,2}"] &\cdots\arrow[l,"\alpha_{N,p_N-2}"] &\circ\arrow[l,"\alpha_{1,p_N-1}"]  & 
      \end{tikzcd}$$
    Note that $\mathrm{End}_{\sC}(T_0)$ is not a hereditary algebra in general, thus the skew-symmetric matrix associated to $Q_{T_0}$ is different from the skew-symmetric matrix $\tilde{B}(\bp,\bla)$ of Euler form.

    (2) If $\X_{\bp,\bla}$ is of parabolic type (see \cite[Section 5.4.1]{Geigle1987}) and each term of $\bp$ is odd (i.e. $\bp=(2r_1+1,2r_2+1)$), then  the cluster category $\sC(\sA_k)$ is triangle equivalent to the cluster category $\sC(\mathrm{mod}kQ)$ of the acyclic quiver $Q$ of type $\tilde{A}_{p_1,p_2}$.  Then $\tilde{B}(\bp,\bla)$ is the same as the skew-symmetric matrix $B_Q$ associated to the quiver $Q$ up to a choice of basis for $\Z^m$. So that  the quantum cluster algebra $\sA(\varLambda,B(\bp,\bla))$ of $\X_{\bp,\bla}$ is isomorphic to the quantum cluster algebra of the acyclic quiver $Q$.

   (3) We define the quantum cluster algebra $\sA(\varLambda,\tilde{B}(\bp,\bla))$ of $\X_{\bp,\bla}$ as a subalgebra of $\hat{\sT}_{\varLambda}$, it follows that each $X_i(t)$ may be a infinite sum of monomials in $\Z[\nu^{\pm}][x_1^{\pm},x_2^{\pm},\cdots, x_m^{\pm}]$. However, we can not deduce that any $X_i(t)$ expressed as a fraction of polynomial of $X_1(t_0), X_2(t_0),\cdots, X_m(t_0)$ is a Laurent polynomial (i.e. the denominator is a monomial). In other words, we do not know whether the quantum cluster algebra $\sA(\varLambda,\tilde{B}(\bp,\bla))$ has the Laurent phenomenon in general. 
   \eremark

 \section{Quantum cluster algebras of \texorpdfstring{$\X_{\bp,\bla}$}{Lg}}\label{sec4}
 \subsection{Cluster multiplication formulas}\label{sec4.1}
   Let $(\varLambda,\tilde{B}(\bp,\bla))$ be a compatible pair. Without loss of generality, assume that $-\varLambda\tilde{B}(\bp,\bla)=I_m$. Let $k=\F_q$ and $v=q^{\frac{1}{2}}$. Recall that $\tilde{\sA}$ is the category $\mathrm{Coh}(\X_{\tilde{\bp},\tilde{\bla}})$ over $k$. 

   To make notations simpler, we will omit the multiplication symbol $*$ of  $H_{\varLambda}(\tilde{\sA})$ and  $\hat{\sT}_{\varLambda,v}$ in the sequel. 

     \begin{lemma}\label{lem2}
       Let $\sM$, $\sN\in \tilde{\sA}$, the following identity holds:
          $$|\mathrm{Gr}_{\underline{e}}(\sM\oplus \sN)|=\sum_{\substack{A,B,C,D,\\ [B]+[D]=\underline{e}}} q^{[B,C]^0}g_{AB}^\sM g_{CD}^{\sN}. $$
    \bproof
     The statement is deduced by applying \cite[Lemma 7]{Hubery2010} to the split exact sequence:  
        $$0\lrw \sN\lrw \sN\oplus \sM \stackrel{\pi}\lrw \sM \lrw 0.$$
    \eproof
     \end{lemma}
    
    Similar to \cite[Theorem 7.4]{Chen2021}, we also have a cluster multiplication formula on the quantum torus $\hat{\sT}_{\varLambda,v}$ specialized at $\nu=v$.
  
    \begin{theorem}\label{thm3.2}
    For $\sM$, $\sN\in \tilde{\sA}$, we have the following equation in  $\hat{\sT}_{\varLambda,v}$:
           \begin{center}
             \begin{align*}
                &(q^{[\sM,\sN]^1}-1)X_{\sM}X_{\sN}=q^{\frac{1}{2}\varLambda(\bm^*,\bn^*)}\sum_{\sL\neq [\sM\oplus \sN]} |\operatorname{Ext}^1_{\tilde{\sA}}(\sM,\sN)_{\sL}| X_{\sL}\\
                &\ \ \ \ \ \ +\sum_{[\sG],[\sF]\neq [\sN]} q^{\frac{1}{2}\varLambda((\bm-\bg)^*,(\bn+\bg)^*)+\frac{1}{2} \langle \bm-g,\bn \rangle }|_{\sF}\Hom_{\tilde{\sA}}(\sN,\tau\sM)_{\tau\sG}| X_{\sG}X_{\sF},
             \end{align*}
           \end{center}
    where $ |\operatorname{Ext}^1_{\tilde{\sA}}(\sM,\sN)_{\sL}| $ means the number of extension classes whose middle term is isomorphic to $\sL$, $|_{\sF}\Hom_{\tilde{\sA}}(\sN,\tau\sM)_{\tau\sG}|$ meas the homomorphism $f\in \Hom_{\tilde{\sA}}(\sN,\tau\sM)$ such that $\Ker f$ isomorphic to $\sF$ and $\Coker f$ isomorphic to $\tau\sG$.
      \begin{proof}
        Since $X_?:H_{\varLambda}(\tilde{\sA})\to \hat{\sT}_{\varLambda,v}$ is an algebra homomorphism, we have
           $$q^{[\sM,\sN]^1}X_{\sM}X_{\sN}=q^{\frac{1}{2}\varLambda(\bm^*,\bn^*)}\sum_{[\sL]\neq[\sM\oplus\sN]} |\operatorname{Ext}^1_{\tilde{\sA}}(\sM,\sN)_{\sL}|X_{\sL}+q^{\frac{1}{2}\varLambda(\bm^*,\bn^*)}X_{\sM\oplus \sN}.$$ 
        On the other hand,   by Lemma \ref{lem2} we have
          $$|\mathrm{Gr}_{\underline{e}}(\sM\oplus \sN)|=\sum_{\substack{A,B,C,D,\\ [B]+[D]=\underline{e}}} q^{[B,C]^0}g_{AB}^\sM g_{CD}^{\sN}.$$
        Hence, 
         \begin{align*}
          &q^{[\sM,\sN]^1}X_{\sM}X_{\sN}-q^{\frac{1}{2}\varLambda(\bm^*,\bn^*)}\sum_{[\sL]\neq[\sM\oplus\sN]} |\operatorname{Ext}^1_{\tilde{\sA}}(\sM,\sN)_{\sL}|X_{\sL}\\
          &=q^{\frac{1}{2}\varLambda(\bm^*,\bn^*)}\sum_{A,B,C,D}q^{- \langle b+d,a+c \rangle }q^{[B,C]^0}g_{AB}^{\sM}g_{CD}^{\sN}X^{-(b+d)^*-^*(a+c)}.
         \end{align*}

         Set $\sigma:=\sum_{\sF,\sG}q^{\frac{1}{2}\varLambda((\bm-\bg)^*,(\bn+\bg)^*)+\frac{1}{2} \langle \bm-\bg,\bn \rangle }|_{\sF}\Hom_{\tilde{\sA}}(\sN,\tau\sM)_{\tau\sG}|X_{\sG}X_{\sF}.$ Note that $|_{\sF}\Hom_{\tilde{\sA}}(\sN,\tau\sM)_{\tau\sG}|=\sum_{S} a_S g_{\sS\sF}^{\sN}g_{\tau\sG,S}^{\tau\sM}$, then 
         $$\sigma=\sum_{\substack{\sF,\sG,\sS\\K,L,X,Y}} q^t a_{\sS}g_{\sS,\sF}^{\sN}g_{\sG,\tau^{-1}\sS}^{\sM}g_{K,L}^{\sG}g_{X,Y}^{\sF}X^{-(l+y)^*- ^*(k+l)}.$$

         where $t=\frac{1}{2}(\varLambda((\bm-\bg)^*,(\bn+\bg)^*)+ \langle \bm-\bg,\bn \rangle - \langle y,x \rangle - \langle l,k \rangle +\varLambda(-l^*-^*k,-y^*-^*x)).$

         Now let us focus on the exponent $t$. Firstly replace the skew-symmetric form $\varLambda$ by $ \langle , \rangle $ as much as possible. Note $(\underline{\dim} \tau^{-1}(S))^*=\tilde{E}'\underline{\dim} \tau^{-1}(S)=-\tilde{E}\underline{\dim} S=-^*s$, then $\tau^{-1}(s)^*=-^*s$. So we have
         \begin{align*}
        2t&=\varLambda(\bm^*,\bn^*)+\varLambda(\tau^{-1}(s)^*,\bg^*)-\varLambda(\bg^*,\bn^*)+ \langle \bm,\bn \rangle - \langle \bg,\bn \rangle +\\
        &\ \  \varLambda((l+k)^*,(y+x)^*)+ \langle l,x \rangle - \langle y,k \rangle - \langle y,x \rangle - \langle l,k \rangle .\\
        &=\varLambda(\bm^*,\bn^*)+\varLambda(\tau^{-1}(s)^*,\bg^*)+ \langle \bm,\bn \rangle +\varLambda(^*s,\bg^*)- \langle \bg,\bff \rangle +\\
        &\ \  \langle l,x \rangle - \langle y,k \rangle - \langle y,x \rangle - \langle l,k \rangle .\\
        &=\varLambda(\bm^*,\bn^*)+ \langle \bm,\bn \rangle - \langle \bg,\bff \rangle + \langle l,x \rangle - \langle y,k \rangle \\
        &\ \  \langle y,x \rangle - \langle l,k \rangle .
         \end{align*}
         
         Secondly, replace $l$ by $\bg-k$ and $x$ by $\bff-y$, then
         \begin{align*}
          2t&=\varLambda(\bm^*,\bn^*)+ \langle \bm,\bn \rangle - \langle \bg,\bff \rangle + \langle \bg-k,\bff-y \rangle - \langle y,k \rangle \\
             &\ \ - \langle y,\bff-y \rangle - \langle \bg-k,k \rangle .\\
             &=\varLambda(\bm^*,\bn^*)+2 \langle \bm-k,\bn-y \rangle - \langle \bm+y-k,\bn-y+k \rangle .
           \end{align*}
        The second equality is induced by $ \langle y,\bn-d \rangle = \langle y,s \rangle =- \langle \tau^{-1}s,y \rangle =- \langle \bm-\bg,y \rangle $  and $ \langle k,\bm-\bg \rangle = \langle k,\tau^{-1}s \rangle =- \langle s,k \rangle =- \langle \bn-\bff,k \rangle $.

        So 
        \begin{align*}
          \sigma=q^{\varLambda(\bm^*,\bn^*)}\sum_{\substack{\sF,\sG,\sS\\K,L,X,Y}}& q^{ \langle \bm-k,\bn-y \rangle -\frac{1}{2} \langle \bm+y-k,\bn-y+k \rangle } \\
          &a_{\sS}g_{\sS\sF}^{\sN}g_{\sG,\tau^{-1}\sS}^{\sM}g_{KL}^{\sG}g_{XY}^{\sF}X^{-(l+y)^*- ^*(k+l)}.
        \end{align*}

        By the associativity of Hall algebra $H_{\varLambda}(\tilde{\sA})$, we have  
            $$\sum_{\sF}g_{\sS\sF}^{\sN}g_{XY}^{\sF}=\sum_{D} g_{\sS X}^Dg_{DY}^\sN \ \ 
            and\ \  \sum_{\sG}g_{\sG,\tau^{-1}\sS}^{\sM}g_{KL}^{\sG}=\sum_{A} g_{L,\tau^{-1}\sS}^{A}g_{KA}^{\sM}.$$
        Then $\sum_{\sF,\sG,\sS} a_{\sS}g_{\sS\sF}^{\sN}g_{\sG,\tau^{-1}\sS}^{\sM}g_{KL}^{\sG}g_{XY}^{\sF}=\sum_{D,A,\sS} a_{\sS}g_{\sS X}^Dg_{DY}^{\sN} g_{L,\tau^{-1}\sS}^{A}g_{KA}^{\sM}$, it follows that
         $$\sum_{\sS,L,X}a_{\sS}g_{\sS X}^D g_{L,\tau^{-1}\sS}^{A}=\sum_{X,L}|_ X\Hom_{\tilde{\sA}}(D,\tau A)_{\tau L}|= |\Hom_{\tilde{\sA}}(D,\tau A)|=q^{[A,D]^1}.$$ and
        $$(l+y)^*+^*(k+x)=(l+y+\tau^{-1}s)^*+^*(k+x+s)=(a+y)^*+^*(d+x).$$
        We can summarize all $\sS$, $L$ and $X$ of $\sigma$ to get
         $$\sigma=q^{\varLambda(\bm^*,\bn^*)}\sum_{A,D,Y,K}q^{ \langle \bm-k,\bn-y \rangle -\frac{1}{2} \langle \bm+y-k,\bn-y+k \rangle } q^{[a,d]^1}g_{DY}^{\sN} g_{KA}^{\sM}X^{-(a+y)^*-^*(d+x)}.$$
        Replace $A,D,Y,K$ by $B,C,D,A$, we have
        $$\sigma=q^{\varLambda(\bm^*,\bn^*)}\sum_{A,B,C,D}q^{ \langle b,c \rangle -\frac{1}{2} \langle b+d,a+c \rangle } q^{[b,c]^1}g_{CD}^{\sN} g_{AB}^{\sM}X^{-(b+d)^*-^*(a+c)}.$$
        which is exactly $X_{\sM\oplus \sN}$.

        Hence, we have 
        $$(q^{[\sM,\sN]^1})X_{\sM}X_{\sN}=-q^{\frac{1}{2}\varLambda(\bm^*,\bn^*)}\sum_{[\sL]\neq[\sM\oplus\sN]} |\operatorname{Ext}^1_{\sA}(\sM,\sN)_{\sL}|X_{\sL}+\sigma.$$
        If $\sF\cong \sN$, then $f\in_{\sF}\Hom_{\tilde{\sA}}(\sN,\tau(\sM))_{\tau(\sG)}$ must be 0. In this case, $\sG\cong\sM$, so $\sigma=\sigma_{[\sF]\neq [\sN]}+X_{\sM}X_{\sN}=:\sigma_2+X_{\sM}X_{\sN}$ and we complete the proof.
      \end{proof}
     \end{theorem}

     \bcorollary\label{cor4.3}
     For an exchange pair $(T_i,T_i^*)$ in $\sC(\tilde{\sA_k})$ such that $Ext^1_{\tilde{\sA}_k}(T_i^*,T_i)\neq 0$, then we have the following identities:
       \begin{equation}\label{eq4.1}
        X_{T_i^*}X_{T_i}=q^{\frac{1}{2}\varLambda(\bm^*,\bn^*)} X_{E} + q^{\frac{1}{2}(\varLambda(\bm^*,\bn^*)-1)}X_{E'},
      \end{equation}
      \begin{equation}\label{eq4.4}
        X_{T_i}X_{T_i^*}=q^{\frac{1}{2}\varLambda(\bn^*,\bm^*)} X_{E} + q^{\frac{1}{2}(\varLambda(\bn^*,\bm^*)+1)}X_{E'},
      \end{equation}
    where $E$ and $E'$ are the middle terms of the exchange triangles respectively, $\bm$ (resp. $\bn$) is the dimension vector of $T_i^*$ (resp. $T_i$).
    \bproof
     By Theorem \ref{thm7.3}, $\dim_k \mathrm{Ext}^1_{\sC(\tilde{\sA}_k)}(T_i,T_i^*)=1$. Note that 
       $$\mathrm{Ext}^1_{\sC_k}(T_i,T_i^*)=\mathrm{Ext}^1_{\tilde{\sA}_k}(T_i,T_i^*)\oplus \mathrm{Ext}^1_{\tilde{\sA}_k}(T_i^*,T_i),$$ 
     and $\mathrm{Ext}^1_{\tilde{\sA}_k}(T_i^*,T_i)\neq 0$, it follows that $\dim_k \mathrm{Ext}^1_{\tilde{\sA}_k}(T_i^*,T_i)=1$, $\dim_k \mathrm{Ext}^1_{\tilde{\sA}_k}(T_i,T_i^*)=0$.  Since $\Hom_k(T_i,\tau T_i^*)\cong \D\mathrm{Ext}^1_{\tilde{\sA}}(T_i^*,T_i)\cong k$, then $\sF:=\Ker f$ (resp. $\sG:=\tau^{-1}\Coker f$) are the same for any nonzero homomorphism $f: T_i\to \tau T_i^*$. Denoted by $\sS$ the image $\mathrm{Im} f$ of $f$. Following from Theorem \ref{thm3.2} we have that
     \begin{equation}\label{eq4.2}
        X_{T_i^*}X_{T_i}=q^{\frac{1}{2}\varLambda(\bm^*,\bn^*)} X_{E} + q^{\frac{1}{2}\varLambda((\bm-\bg)^*,(\bn+\bg)^*)+\frac{1}{2} \langle \bm-g,\bn \rangle }X_{\sG}X_{\sF},
    \end{equation}
    where $g=\underline{\dim}\sG$ and $f=\underline{\dim}\sF$.  Note that $E'\cong \sF\oplus \sG$ is rigid and $X_?: H_{\varLambda}(\tilde{\sA})\to \sT_{\varLambda,v}$ is an algebra homomorphism, $X_{\sG}X_{\sF}=q^{\frac{1}{2}\varLambda(g^*,f^*)}X_{E'}$. Comparing Equation (\ref{eq4.2}) with Equation (\ref{eq4.1}), it suffices to show that
       $$\varLambda((\bm-\bg)^*,(\bn+\bg)^*)+ \langle \bm-g,\bn \rangle+\varLambda(g^*,f^*)=\varLambda(\bm^*,\bn^*)-1.$$
    Using $\bm=\tau^{-1}s+g$, $\bn=f+s$ and $\tau^{-1}(s)^*=-^*s$, where $s=\underline{\dim}\sS$, we have that 
       \begin{align*}
        &\varLambda((\bm-\bg)^*,(\bn+\bg)^*)+ \langle \bm-g,\bn \rangle+\varLambda(g^*,f^*) \\           
        &=\varLambda(\bm^*,\bn^*)+\varLambda(\tau^{-1}(s)^*,g^*)-\varLambda(g^*,\bn^*)+ \langle \bm,\bn \rangle- \langle g,\bn \rangle+\varLambda(g^*,f^*),\\
        &=\varLambda(\bm^*,\bn^*)+\varLambda(\tau^{-1}(s)^*,g^*)+\varLambda(^*s,g^*)+ \langle \bm,\bn \rangle- \langle g,f \rangle,\\
        &=\varLambda(\bm^*,\bn^*)+ \langle \bm,\bn \rangle- \langle g,f \rangle.
         \end{align*}
      Applying $\Hom_{\tilde{\sA}_k}(\sF,-)$ to the exact sequence $\sF\rightarrowtail T_i\twoheadrightarrow \sS$ and note $\mathrm{Ext}^1_{\tilde{\sA}_k}(\sF,T_i)=0$, we obtain $\D\Hom_{\tilde{\sA}}(\sS,\tau\sF)\cong$ $\mathrm{Ext}^1(\sF,\sS)=0$.  Then apply $\Hom_{\tilde{\sA}_k}(-,\sF)$ to the exact sequence $\tau^{-1}S\rightarrowtail T_i^*\twoheadrightarrow \sG$, we deduce that $\Hom_{\tilde{\sA}_k}(\sG,\sF)\cong \Hom_{\tilde{\sA}_k}(T_i^*,\sF)$. Finally apply $\Hom_{\tilde{\sA}_k}(T_i^*,-)$ to the exact sequence $\sF\rightarrowtail T_i\twoheadrightarrow \sS$  to get $\Hom_{\tilde{\sA}_k}(T_i^*,\sF)\cong \Hom_{\tilde{\sA}_k}(T_i^*,T_i)$ since  $\Hom_{\tilde{\sA}_k}(T_i^*,\sS)\rightarrowtail \Hom_{\tilde{\sA}_k}(T_i^*,\tau T_i^*)\cong D\mathrm{Ext}_{\tilde{\sA}_k}^1(T_i^*,T_i^*)=0$. Thus $\langle g,f\rangle=\dim_k \Hom_{\tilde{\sA}_k}(\sG,\sF)=\langle \bm,\bn\rangle +1$, implying that
      $$\varLambda(\bm^*,\bn^*)+ \langle \bm,\bn \rangle- \langle g,f \rangle=\varLambda(\bm^*,\bn^*)-1,$$
      which gives rise to the first equation.

     For the second equation, we have $X_{T_i}X_{T_i^*}=q^{\frac{1}{2}\varLambda(\bn^*,\bm^*)}X_{T_i\oplus T_i^*}$ for $\mathrm{Ext}^1_{\tilde{\sA}_k}(T_i,T_i^*)=0$. On the other hand, from $\mathrm{Ext}^1_{\tilde{\sA}_k}(T_i^*,T_i)\cong k$ we have that
        $$qX_{T_i^*}X_{T_i}=q^{\frac{1}{2}\varLambda(\bm^*,\bn^*)}(q-1)X_E+q^{\frac{1}{2}\varLambda(\bm^*,\bn^*)} X_{T_i\oplus T_i^*}.$$
    Combining with Equation $(\ref{eq4.1})$ and $\varLambda$ is skew-symmetric, we will obtain Equation $(\ref{eq4.4})$.
    \eproof

     \ecorollary

     \bexample\label{ex3.3}
      Take $\X_{\bp,\bla}$ to be the projective line $\P^1$, then $\tau(\sF)=\sF(-2)$. The matrix $E$ of the Euler form on $K_0(\sA)$ under the basis $\{\hat{\sO},\hat{S}_x\}$  is $E=\begin{bmatrix}1  &1\\ -1  &0\end{bmatrix}$ and $\varLambda'=\begin{bmatrix}0  &-1/2\\ 1/2  &0\end{bmatrix}$. Then $-2\varLambda' B=I_2$. $\varLambda=2\varLambda'$ and $\tilde{\sA}=\sA$ is $\mathrm{Coh}(\P^1)$ over $\F:=\F_{q^2}$.

      (1) Take $\sM=\sO(2)$, $\sN=\sO$. Their dimension vectors are $m=\begin{bmatrix}1  \\2 \end{bmatrix}$ and $n=\begin{bmatrix}1  \\0 \end{bmatrix}$ respectively. Then 
       \begin{center}
       $ m^*=\begin{bmatrix}-1  \\1 \end{bmatrix}$,   $n^*=\begin{bmatrix}1  \\1 \end{bmatrix}$,  and $\varLambda(m^*,n^*)=2$. 
        \end{center}
      The only non-trivial extension of $\sO$ by $\sO(2)$ is 
        $$0\to \sO\to \sO(1)^{\oplus 2}\to \sO(2)\to 0.$$
      Note for any nonzero homomorphism $f\in \Hom_{\sA}(\sO,\tau(\sO(2)))\cong \F_{q^2}$, $f$ is isomorphic, thus Theorem \ref{thm3.2} applied to $\sM,\sN$ is 
          \begin{equation}\label{eqn3.1}
            (q^2-1)X_{\sO(2)}X_{\sO}=q(q^2-1)X_{\sO(1)^{\oplus 2}}+(q^2-1).
          \end{equation}
      On the other hand, the cluster character of the vector bundle $\sO(l)$ is
          $$X_{\sO(l)}=X^{-(\begin{smallmatrix}l+1\\ -1\end{smallmatrix})}+\sum_{r\geq l}q^{-(l+r)}[l+r+1]_{q^2}X^{-(\begin{smallmatrix}l+2r+1\\ 1\end{smallmatrix})}.$$
      By direct computation, we have the following identities in $\hat{\sT}_{\varLambda,v}$.
          \begin{equation}\label{eqn3.2}
            X_{\sO(2)}X_{\sO}=qX_{\sO(1)}X_{\sO(1)}+1.
          \end{equation}
      Note $[\sO(1)]*[\sO(1)]=[\sO(1)^{\oplus 2}]$ in $H_{\varLambda}(\sA)_{\F}$, implying $X_{\sO(1)}^2=X_{\sO(1)^{\oplus 2}}$. So from identity (\ref{eqn3.2}), we have
      \begin{equation}
        X_{\sO(2)}X_{\sO}=qX_{\sO(1)^{\oplus 2}}+1.
      \end{equation}
      which gives the identity (\ref{eqn3.1}).

      (2) Take $\sM=S_x$, $\sN=\sO$, where $\mathrm{deg}(x)=1$. Then
        \begin{center}
           $ m^*=\begin{bmatrix}-1  \\0 \end{bmatrix}$, and  $n^*=\begin{bmatrix}1  \\1 \end{bmatrix}$,  $\varLambda(m^*,n^*)=1$. 
         \end{center}
      Notice $|\operatorname{Ext}^1_{\sA}(S_x,\sO)|=q^2-1$ and any nonzero homomorphism $g\in \Hom_{\sA}(\sO,\tau(S_x))\cong k$ is surjective with $\Ker g=\sO(-1)$. Hence the quantum cluster multiplication formula given in Theorem \ref{thm3.2} applied to $\sM,\sN$ is 
         \begin{equation}
           (q^2-1)X_{S_x}X_{\sO}=q^{\frac{1}{2}}(q^2-1)X_{\sO(1)}+q^{-\frac{1}{2}}(q^2-1)X_{\sO(-1)}.
         \end{equation}
      On the other hand, the cluster character of $S_x$ is 
        $$X_{S_x}=X^{-(\begin{smallmatrix}1\\ 0\end{smallmatrix})}+X^{-(\begin{smallmatrix}-1\\ 0\end{smallmatrix})}.$$
      By direct computations, we have
         \begin{equation}\label{eqn3.4}
           X_{S_x}X_{\sO(m)}=q^{\frac{1}{2}}X_{\sO(m+1)}+q^{-\frac{1}{2}}X_{\sO(m-1)}.
         \end{equation}
      which gives rise to the equation (\ref{eqn3.4}).
     \eexample
    
  \subsection{Quantum F-polynomials}\label{sec4.2} In this subsection, we still assume that $\varLambda \tilde{B}(\bp,\bla)=-I_m$ for some skew-symmetric matrix $\varLambda$ of integers. Write $\tilde{B}:=\tilde{B}(\bp,\bla)$. Let $k$ be a finite field. Recall that we have constructed a valued regular $m$-tree $\sT_m$ in Section \ref{sec3.4}. Remind that $T(t)$ is a symbol associated to $t\in \T_m$ such that $T(t)^k$ is the cluster-tilting object in $\sC(\tilde{\sA}_k)$. In this subsection, we want to show that $\mathrm{Gr}_{e}(T_i(t)^k)$ is a polynomial of $|k|$ for $t\in \T_m$, $1\leq i\leq m$ and $\bee\in \Z^m$.

    The quantum cluster character of $\sF^k\in \tilde{\sA}_k$ is 
        $$X_{\sF^k}=\sum_{\bee\leq \bff} q^{-\frac{1}{2}\langle \bff-\bee,\bee \rangle}|\mathrm{Gr}_{\bee}(\sF^k)|_q X^{-\bff^*+\bee^*-^*\bee},$$
    where $\bff=\underline{\dim}\sF^k$.  Recall that we have defined a partial order associated to $\bff$ on $\{-^*\bff+\tilde{B}\bee|\ \bee\leq \bff\}$ in Section \ref{sec3.4}.   One observation is that each $X_{T_j(t)^k}$ has a unique maximal degree for $1\leq j\leq m$ and $t\in \T_m$.

    \btheorem
      For $T_i(t')$ with $t'\in \T_m$, $1\leq i\leq m$, there exists a  $\Z$-polynomial $P(z)$ such that the cardinality $|\mathrm{Gr}_{\bee}(T_i(t')^k)|=P(|k|^{\frac{1}{2}})$.
     \bproof
      We prove the statement by induction on $t$ from the root $t_0\in \T_m$. Notice that we have already shown that $|\mathrm{Gr}_{\bee}(T_i(t_0)^k)|$ is a polynomial of $|k|^{\frac{1}{2}}$ in Section \ref{sec3.4}.  Assume that for $\begin{tikzcd}[column sep=huge]t\arrow[r, no head,"{(a_{ij}(t),b_{ij}(t))}"] &t'\end{tikzcd}$, the statement holds for $T_j(t)$ for any $1\leq j\leq m$. Namely, $|\mathrm{Gr}_{\bee}(T_j(t)^k)|$ is a $\Z$-polynomial of $|k|^{\frac{1}{2}}$ for any $\bee \leq d_j(t)$. Using Corollary \ref{cor4.3}, say one of equations is
    \begin{equation*}
       X_{T_i(t')^k}X_{T_i(t)^{k}}=q^{\frac{1}{2}\varLambda(d_i(t')^*,d_i(t)^*)}q^r\prod_{j\neq i} X_{T_j(t)^k}^{a_{ij}} + q^{\frac{1}{2}(\varLambda(d_i(t')^*,d_i(t)^*)-1)}q^{r'}\prod_{j\neq i} X_{T_j(t)^k}^{b_{ij}},
    \end{equation*}
    comparing degrees from the unique maximal one on both side, we can calculate the cardinality $|\mathrm{Gr}_{\boldsymbol{e}}(T_i(t')^{k})|$ for each $\boldsymbol{e}\leq d_i(t')$. Note that $|\mathrm{Gr}_{\bee'}(T_i(t)^k)|$ and $|\mathrm{Gr}_{\bee'}(T_j(t)^{k})|$ for $j\neq i$ are $\Z$-polynomials of $|k|$ by induction hypothesis, it follows that  $|\mathrm{Gr}_{\bee}(T_i(t')^k)|$ is equal to ${v^{-s}}f(v)\in \Z[v^{\pm1}]$ where $v=|k|^{\frac{1}{2}}$. Because $|\mathrm{Gr}_{\bee}(T_i(t')^k)|$ is an integer for any $|k|>2$, we have $|\mathrm{Gr}_{\bee}(T_i(t')^k)|\in \Z[v]$. Otherwise $t^{-s}f(t)=f_1(t)+f_2(t^{-1})$ with $f_1\in \Z[t]$ and $f_2\in t^{-1}\Z[t^{-1}]$. when $t=q^r$ goes to $+\infty$, $f_1(t)\in \Z$ while $|f_2(t^{-1})|\leq 1$. This is contradict to $f(t)\in \Z$. The proof is completed.
     \eproof
    \etheorem

    Recall the quantum cluster character of $T_i(t)^k$ is
        $$   X_{T_i(t)^k}=\sum_{\bee\leq d_i(t)} |k|^{-\frac{1}{2} \langle d_i(t)-\bee,\bee \rangle }|\mathrm{Gr}_{\bee}(T_i(t)^k)| X^{-(d_i(t)-\bee)^*-^*\bee}.$$
    From the last theorem, we know that $|\mathrm{Gr}_{\bee}(T_i(t)^k)|$ is a $\Z$-polynomial of $|k|^{\frac{1}{2}}$, it follows that there exists a unique element $X_{T_i(t)} \in \hat{\sT}_{\varLambda}$ such that
        $$X_{T_i(t)}|_{\nu=|k|^{\frac{1}{2}}}= X_{T_i(t)^k}.$$
    By Corollary \ref{cor4.3}, for an exchange pair $(T_i(t)^k,T_i(t')^k)$ with $\begin{tikzcd}[column sep=huge]t\arrow[r, no head,"{(a_{ij}(t),b_{ij}(t))}"] &t'\end{tikzcd}$, if $\mathrm{Ext}^1_{\tilde{\sA}_k}(T_i(t)^k,T_i(t')^k)=0$, then we have that
    \begin{equation}\label{eq4.10}
      X_{T_i(t')^k}X_{T_i(t)^k}=q^{\frac{1}{2}\varLambda(d_i(t')^*,d_i(t)^*)} q^s \prod_{j\neq i} X_{T_j(t)^k}^{a_{ij}} + q^{\frac{1}{2}(\varLambda(d_i(t')^*,d_i(t)^*)-1)}q^{s'} \prod_{j\neq i} X_{T_j(t)^k}^{b_{ij}}.
    \end{equation}
    by observing that $E^k=\bigoplus_{j\neq i} T_j(t)^{k,\oplus a_{ij}}$ and $E^{'k}=\bigoplus_{j\neq i} T_j(t)^{k,\oplus b_{ij}}$ are rigid,  where $s$ and $s'$ are the same as  Definition \ref{def3.10}. If $\mathrm{Ext}^1_{\tilde{\sA}_k}(T_i(t')^k,T_i(t)^k)=0$, then we have that
    \begin{equation}\label{eq4.11}
      X_{T_i(t')^k}X_{T_i(t)^k}=q^{\frac{1}{2}\varLambda(d_i(t')^*,d_i(t)^*)+1} q^s \prod_{j\neq i} X_{T_j(t)^k}^{a_{ij}} + q^{\frac{1}{2}(\varLambda(d_i(t')^*,d_i(t)^*))}q^{s'} \prod_{j\neq i} X_{T_j(t)^k}^{b_{ij}}.
    \end{equation}
    
    Moreover, $T(t)^k=\bigoplus_{1\leq j\leq m} T_j(t)^k$ is rigid, it follows that 
       \begin{equation}\label{eq4.12}
        X_{T_j(t)^k}X_{T_i(t)^k}=q^{\varLambda(d_j(t)^*,d_i(t)^*)} X_{T_i(t)^k}X_{T_j(t)^k}.
       \end{equation}
       
   Let $\T_n(\bp,\bla)$ be the subgraph of $\T_m$ associated to $(\bp,\bla)$ given in Definition \ref{def3.6}.  Note that $X_{T_i(t_0)}=X_i(t_0)$, using Equations (\ref{eq4.10}) and $(\ref{eq4.11})$ we reach the following
    \btheorem\label{thm4.6}
     The quantum cluster algebra $\sA(\varLambda,\tilde{B}(\bp,\bla))$ is a $\Z[\nu^{\pm1}]$-subalgebra of $\hat{\sT}_{\varLambda}$ generated by $X_{T_i(t)}$ for $t\in \T_n(\bp,\bla)$, $1\leq i\leq n$ and $X_{T_l(t_0)}^{\pm}$ for $n< l\leq m$. 
    \etheorem

    There is a bar-involution $\overline{\bullet}$ on the complete quantum torus $\hat{\sT}_{\varLambda}$, given by $\nu^{\pm} \mapsto \nu^{\mp}$, and $X^{\alpha} \mapsto X^{\alpha}$, $\alpha\in \Z^m$. It is clear that $X_{T_i(t_0)}$ and $X_{T_l(t_0)}^{\pm}$ is bar-invariant (i.e. $\overline{X_{T_i(t_0)}}$=$X_{T_i(t_0)}$)  for each $1\leq i\leq n$, $n<l\leq m$ .

    \bproposition{prop4.7}
    For $1\leq i\leq n$ and $t\in \T_n(\bp,\bla)$, $X_{T_i(t)}$  is bar-invariant.
    \bproof
    For $\begin{tikzcd}[column sep=huge]t_0\arrow[r, no head,"{(a_{ij}(t_0),b_{ij}(t_0))}"] &t\end{tikzcd}$, if $\mathrm{Ext}^1(T_i(t_0),T_i(t))=0$, then by Theorem \ref{thm4.6}  and Corollary \ref{cor4.3} we have
    \begin{equation}
      X_{T_i(t)}X_{T_i(t_0)}=\nu^{\varLambda(d_i(t)^*,d_i(t_0)^*)} X_{E} + \nu^{(\varLambda(d_i(t)^*,d_i(t_0)^*)-1)}X_{E'}.
    \end{equation}
   Note that $E$ and $E'$ are rigid and $X_{T_j(t_0)}$ is bar-invariant, $X_E$ and $X_{E'}$ are bar-invariant. Applying $\overline{\bullet}$ to the equation above,
    \begin{equation*}
      X_{T_i(t_0)}\overline{X_{T_i(t)}}=\nu^{\varLambda(d_i(t_0)^*,d_i(t)^*)} X_{E} + \nu^{(\varLambda(d_i(t_0)^*,d_i(t)^*)+1)}X_{E'}.
    \end{equation*}
    On the other hand, using Corollary \ref{cor4.3} again, we also have
    \begin{equation*}
      X_{T_i(t_0)}X_{T_i(t)}=\nu^{\varLambda(d_i(t_0)^*,d_i(t)^*)} X_{E} + \nu^{(\varLambda(d_i(t_0)^*,d_i(t)^*)+1)}X_{E'},
    \end{equation*}
    Therefore, $X_{T_i(t_0)}\overline{X_{T_i(t)}}=X_{T_i(t_0)}X_{T_i(t)}$ in $\hat{\sT}_{\varLambda}$, it follows that $\overline{X_{T_i(t)}}=X_{T_i(t)}$. 

    Repeat last procedure, we can prove that $X_{T_i(t)}$ is bar-invariant for any $t\in \T_n(\bp,\bla)$.
    \eproof
    \eproposition

    \bcorollary\label{co4.8}
    The map 
    $$\begin{aligned}
      \overline{\bullet}: \sA(\varLambda,\tilde{B}(\bp,\bla)) &\lrw \sA(\varLambda,\tilde{B}(\bp,\bla)),\\
      \nu^{\pm} &\mapsto \nu^{\mp},\\
      X_{T_i(t)}&\mapsto X_{T_i(t)}.
    \end{aligned}$$ 
    give rise to a bar-involution.
    \ecorollary

  \subsection{Specialized quantum cluster algebras}\label{sec4.3}
       Assume  $-\varLambda B(\tilde{\bp},\tilde{\bla})=dI_m$ where $\tilde{B}:=B(\tilde{\bp},\tilde{\bla})$ is invertible. Namely, each $\tilde{p}_i$ is odd. Let $k=\F_{q^d}$ and $v=q^{\frac{1}{2}}$. Recall that the specialized quantum cluster algebra $\sA_q(\varLambda,\tilde{B})$ of $\X_{\bp,\bla}$ at $\nu=v$ is the subalgebra of $\hat{\sT}_{\varLambda,v}$ generated by quantum cluster characters $X_{T_i(t)^k}$ of $T_i(t)^k$ for $t\in \T_n(\bp,\bla)$, $1\leq i\leq m$ and $X_{T_j(t_0)^k}^{\pm}$ for $n<j\leq m$
       
       Let $\mathrm{CH}'_{\varLambda}(\tilde{\sA}_k)$ be the subalgebra of $H_{\varLambda}(\tilde{\sA}_k)$ generated by $[\sO(l\vec{c})^k]$, and $[S_{i,j}^k]$,  for $l\in\Z $, $1\leq i\leq N$, and $1\leq j\leq p_i$. Consider $\mathrm{CH}'_{\varLambda}(\tilde{\sA}_k)\otimes_{\Z[v^{\pm1}]} \Q$, note that we have the following exact sequences
          $$0\lrw \sO((j-1)\vec{x}_i)^k\lrw \sO(j\vec{x}_i)^k\lrw S_{i,j}^k\lrw 0,$$ 
      it follows that $\sO(\vec{l})^k\in \mathrm{CH}'_{\varLambda}(\tilde{\sA}_k)\otimes_{\Z[v^{\pm1}]} \Q$ for any $\vec{l}\in \L(\tilde{\bp},\tilde{\bla})$.
     
       An object $\sF$ in a hereditary category $\sA$ is called exceptional if $\sF$ is rigid and $\mathrm{End}_{\sA}(\sF)$ is a division ring. In addition, it is well known that  $\mathrm{End}_{\tilde{\sA}_k}(\sF)\cong k$ for an indecomposable rigid object $\sF\in \tilde{\sA}_k$ (see \cite[Proposition 6.4.2]{Chen2009}). Thus exceptional objects in $\tilde{\sA}_k$ are precisely indecomposable rigid objects in $\tilde{\sA}_k$.  So the following theorem will be applied to finite fields \cite{CrawleyBoevey1992,Lenzing2011, Meltzer1995,Lenzing2009}.
      
       \btheorem[{\cite[Theorem 1]{Kedzierski2013}}]\label{thm4.7}
       Let $\sF\in \sA$ be an exceptional vector bundle of rank greater than one on a weighted projective line $\X_{\bp,\bla}$ over an algebraically closed field. Then there are exceptional objects $\sF'$ and $\sF''$ with the following properties:
       \begin{itemize}
        \item [(i)] $\Hom_{\sA}(\sF',\sF'') = \Hom_{\sA}(\sF'',\sF') = \mathrm{Ext}_{\sA}^1(\sF',\sF'')=0$, and there is a nonsplit exact sequence
            $$0\lrw \sF^{'\oplus a} \lrw \sF \lrw  \sF^{''\oplus b} \lrw 0, $$
        where $(a,b)$ is the dimension vector of $\sF\in \sC(\sF',\sF'')\simeq  \mathrm{mod}k\Theta_r$ and $r=\dim \mathrm{Ext}^1_{\sA}(\sF'',\sF')$. 
        \item [(ii)] $\mathrm{rank}(\sF')< \mathrm{rank}(\sF)$ and $\mathrm{rank}(\sF'')< \mathrm{rank}(\sF)$.
       \end{itemize}
       Here  $\sC(\sF',\sF'')\subset \sA$ is the subcategory containing $\sF$ and $\sF''$, closed under extensions, kernels of epimorphisms, and cokernels of monomorphisms. $\Theta_r$ denotes the $r$-Kronecker quiver:
         $$\begin{tikzcd}
          1  &2\arrow[l,shift right=4,"\vdots"]\arrow[l, shift right=2]\arrow[l, shift left=3] \ \ (r \text{ arrows}).
         \end{tikzcd}$$
       \etheorem

     \begin{corollary}\label{cor4.8}
      Each indecomposable rigid object in $\tilde{\sA}_k$ belongs to $\mathrm{CH}'_{\varLambda}(\tilde{\sA}_k)\otimes_{\Z[v^{\pm1}]} \Q$.
      \bproof
      Indecomposable rigid objects in $\tilde{\sA}_k$ are exceptional vector bundles and exceptional torsion sheaves which lie in $\mathrm{Tor}_{\lambda_i}$ defined in Lemma \ref{lem2.1}. It is clear that all exceptional torsion sheaves and line bundles belong to $\mathrm{CH}'_{\varLambda}(\tilde{\sA}_k)\otimes_{\Z[v^{\pm1}]} \Q$. For an exceptional vector bundle $\sF^k$, by Theorem \ref{thm4.7}, it suffices to show that $\sF^k$ of rank $\geq 2$ belongs to  the $\Q$-linear composition Hall algebra $\mathrm{CH}_{\varLambda}(\sC(\sF',\sF''))\otimes_{\Z[v^{\pm1}]} \Q$. Note that $\sC(\sF',\sF'')$ in Theorem \ref{thm4.7} is equivalent to the module category $\mathrm{mod}k\Theta_r$, it is equivalent to show that for each indecomposable rigid object $M$ in $\mathrm{mod}k\Theta_r$, $[M]\in \mathrm{CH}(\mathrm{mod}k\Theta_r)\otimes_{\Z[v^{\pm1}]} \Q$. Without loss of generality, we assume $r=2$, namely, $\Theta_r$ is the Kronecker quiver. Now the statement is induced from \cite[Theorem 1]{Zhang1996}. Since every indecomposable object in $\mathrm{mod}k\Theta_2$ is either an indecomposable preprojective or an indecomposable preinjective object.
      \eproof
     \end{corollary}
    
     Let $I_k$ be a 2-sided ideal of $\mathrm{CH}'_{\varLambda}(\tilde{\sA}_k)\otimes_{\Z[v^{\pm1}]} \Q$ generated by 
      \begin{equation}\label{eq4.13}
         [T_j(t)^k][T_i(t)^k]-v^{2\varLambda(d_j(t)^*,d_i(t)^*)}[T_i(t)^k][T_j(t)^k], \text{ if } j\neq i,
      \end{equation}
      and if $\begin{tikzcd}[column sep=huge]t\arrow[r, no head,"{(a_{ij}(t),b_{ij}(t))}"] &t'\end{tikzcd}$ in $\T_m$ such that $\mathrm{Ext}^1_{\tilde{\sA}_k}(T_i(t)^k,T_i(t')^k)=0$,
      \begin{equation}\label{eq4.14}
      [T_i(t')^k][T_i(t)^k]-v^{\varLambda(d_i(t')^*,d_i(t)^*)} v^s \prod_{j\neq i}[T_j(t)^k]^{a_{ij}(t)}  
      - v^{\varLambda(d_i(t')^*,d_i(t)^*)-1} v^{s'} \prod_{j\neq i}[T_j(t)^k]^{b_{ij}(t)},
      \end{equation}  
      or if  $\begin{tikzcd}[column sep=huge]t\arrow[r, no head,"{(a_{ij}(t),b_{ij}(t))}"] &t'\end{tikzcd}$ in $\T_m$ such that $\mathrm{Ext}^1_{\tilde{\sA}_k}(T_i(t')^k,T_i(t)^k)=0$,
      \begin{equation}\label{eq4.15}
        [T_i(t')^k][T_i(t)^k]-v^{\varLambda(d_i(t')^*,d_i(t)^*)+1} v^s \prod_{j\neq i}[T_j(t)^k]^{a_{ij}(t)}  
        - v^{\varLambda(d_i(t')^*,d_i(t)^*)} v^{s'} \prod_{j\neq i}[T_j(t)^k]^{b_{ij}(t)},
        \end{equation}  
    
     For an element $[\sF^k]$ in $\mathrm{CH}'_{\varLambda}(\tilde{\sA}_k)\otimes_{\Z[v^{\pm1}]} \Q$, we still denote $[\sF^k]$ by the image of $[\sF^k]$ in the quotient algebra $(\mathrm{CH}'_{\varLambda}(\tilde{\sA}_k)\otimes_{\Z[v^{\pm1}]} \Q)/I_k$.

      \btheorem\label{thm4.9}
      There is a homomorphism of algebras :
        $$\phi_k: (\mathrm{CH}'_{\varLambda}(\tilde{\sA}_k)\otimes_{\Z[v^{\pm1}]} \Q)/I_k  \lrw \sA_q(\varLambda,B(\tilde{\bp},\tilde{\bla}))\otimes_{\Z[v^{\pm1}]} \Q,$$
      which maps $[T_i(t)^k]$ to $X_{T_i(t)^k}$ for $1\leq i\leq m$ and $t\in \T_m$. In particular, $\phi_k$ is an isomorphism.
      \bproof
      By Equations (\ref{eq4.10})-(\ref{eq4.12}) and $X_?: H_{\varLambda}(\tilde{\sA}_k)\to \sA_q(\varLambda,\tilde{B})$ is an algebra homomorphism, it can be seen that $\phi_k$ is a homomorphism of algebras. Moreover, following from Corollary \ref{cor4.8}, $\phi_k$ is surjective. On the other hand, the defining relations in $\sA_q(\varLambda,\tilde{B})$ are exactly relations (\ref{eq4.13})-(\ref{eq4.15}) by Lemma \ref{lem3.9}, which induces a homomorphism 
        $$\psi_k:\sA_q(\varLambda,\tilde{B})\otimes_{\Z[v^{\pm1}]} \Q\lrw (\mathrm{CH}'_{\varLambda}(\tilde{\sA}_k)\otimes_{\Z[v^{\pm1}]} \Q)/I_k,$$
      mapping $X_{T_i(t)^k}$ to $[T_i(t)^k]$. As a consequence, we have $\phi_k\psi_k=\id$ and $\psi_k\phi_k=\id$, which means that $\phi_k$ is an isomorphism.
      \eproof
      \etheorem

  \section{The quantum cluster algebra of \texorpdfstring{$\P^1$}{Lg}}\label{sec5}
     In this section, we study bases of the quantum cluster algebra $\sA(\varLambda,B)$ of $\P^1$

  \subsection{ \texorpdfstring{$\sA(\varLambda,B)$}{Lg} and \texorpdfstring{$\sA(2,2)$}{Lg}}\label{sec5.1}\ \ The compatible pair $(\varLambda,B)$ associated to $\P^1$ is 
    $$\varLambda=\begin{bmatrix}0  &-1\\ 1  &0\end{bmatrix},\ \ \ B=\begin{bmatrix}1  &-1\\ 1  &0\end{bmatrix}-\begin{bmatrix}1  &1\\ -1  &0\end{bmatrix},\ \ \text{and}\ \ -\varLambda B=2I_2.$$
   Note that all indecomposable rigid objects in $\mathrm{Coh}(\P^1)$ are line bundles and exceptional pairs are exactly $(\sO(l),\sO(l+1))$ and $(\sO(l+1),\sO(l))$. Then $\sA(\varLambda,B)$ is the subalgebra of $\hat{\sT}_{\varLambda}$ generated by $X_{\sO(l)}$ for $l\in \Z$, subject to 
     \begin{align}
       \label{id6.1}  X_{\sO(l+2)}X_{\sO(l)}&=\nu^2 X_{\sO(l+1)}^2+1,\\
        \label{id6.2} X_{\sO(l)}X_{\sO(l+2)}&=\nu^{-2} X_{\sO(l+1)}^2+1,\\
       \label{id6.3} X_{\sO(l+1)}X_{\sO(l)}&=\nu^2X_{\sO(l)}X_{\sO(l+1)}.
       \end{align}
    where
    $$X_{\sO(l)}=X^{-(\begin{smallmatrix}l+1\\ -1\end{smallmatrix})}+\sum_{r\geq l}\nu^{-2(l+r)}[l+r+1]_{\nu^4}X^{-(\begin{smallmatrix}l+2r+1\\ 1\end{smallmatrix})}.$$ 
    Note that we have defined a bar-involution $\overline{\bullet}$ on $\sA(\varLambda,B)$ in Corollary \ref{co4.8}, the Identity (\ref{id6.2}) can be induced by applying $\overline{\bullet}$ to Identity (\ref{id6.1}). Thus, we will omit Identity (\ref{id6.2}) as defining relations in the sequel.

    Recall that the compatible pair $(\varLambda(\Theta_2),B(\Theta_2))$ of the Kronecker quiver $\Theta_2: \begin{tikzcd} 1 &2\arrow[l,shift left]\arrow[l,shift right]  \end{tikzcd}$ is
    $$B(\Theta_2)=\begin{bmatrix}0  &-2\\ 2  &0\end{bmatrix}\ \ \ \text{and}\ \ \ \varLambda(\Theta_2)=\begin{bmatrix}0  &-1\\ 1  &0\end{bmatrix},$$
  which is the same as $(\varLambda,B)$ of $\P^1$. The quantum cluster algebra $\sA(2,2)$ of the Kronecker quiver defined in \cite{Ding2012b} is the  subalgebra of $\sT_{\varLambda}$ generated by $X_{V(l)}$ for $l\in \Z$ subject to 
      $$\begin{aligned}
      X_{V(l-1)}X_{V(l)}&=\nu^2X_{V(l)}X_{V(l-1)},\\
      X_{V(l-2)}X_{V(l)}&=\nu^2 X_{V(l-1)}^2+1,
      \end{aligned}$$
  where $V(l)$ is the indecomposable preprojective $k\Theta_2$-module $P_l=(1-l,-l)$ for $l\leq  0$, the indecomposable preinjective $k\Theta_2$-module $I_{l-2}=(l-3,l-2)$ for $l\geq  -$, and $V(1)=P_2[1]$, $V(2)=P_1[1]$ for $l=1,2$ by setting $X_{P_l[1]}=x_l$. By the definition of $\sA(2,2)$, we have the following

    \bproposition
     There is an isomorphism of algebras :
          $$    \kappa :\sA(2,2)\lrw  \sA(\varLambda,B), \qquad X_{V(l)}\mapsto X_{\sO(-l)}.$$
    \eproposition
    \bremark
    In Remark \ref{rem2}, we state that the quantum cluster algebra of $\X_{\bp,\bla}$ is isomorphic to the quantum cluster algebra of acyclic quiver $Q$ of type $\tilde{A}_{p_1,p_2}$ when $\bp=(2r_1+1,2r_2+1)$. In particular, in the case of $\bp=(1,1)$, $\X_{\bp,\bla}$ is $\P^1$ and $Q$ is the Kronecker quiver. The isomorphism from $\sA(2,2)$ to $\sA(\varLambda,B)$ is given in the last Proposition. In general, we can explicitly give the isomorphism between them through the equivalence $\sC(\sA_k)\simeq \sC(\mathrm{mod}kQ)$ of cluster categories.
    \eremark
  
  Denote $X_{n\delta}:=X_{S_{x_n}}=X^{-(\begin{smallmatrix}n\\ 0\end{smallmatrix})}+X^{-(\begin{smallmatrix}-n\\ 0\end{smallmatrix})}$ for some simple torsion sheaf $S_{x_n}$ supported on the point $x_n\in \P^1$ of degree n. In the sequel, we will show $X_{n\delta}\in \sA(\varLambda,B)$. We will call the subalgebra of $\hat{\sT}_{\varLambda}$ generated by $X_{n\delta}$, $n\in \N$ the $torsion$-$part$, and denote it by $\sA^{tor}(\varLambda,B)$. Obviously, if we show $X_{n\delta}\in \sA(\varLambda,B)$, then  $\sA^{tor}(\varLambda,B)\subset \sA(\varLambda,B)$.
  
\subsection{Bases of the torsion part}
   Let $E_x^{(n)}\in \sA_{\F}$ be the indecomposable torsion sheaf of length $n$ supported on $x\in \P$ with $\mathrm{deg}(x)=1$.  Then, the quantum cluster character of $E_x^{(n)}$ is 
      $$X_{E_x^{(n)}}=\sum_{l=0}^n X^{-(\begin{smallmatrix}n-2l \\ 0 \end{smallmatrix})}.$$
  Note that the quantum cluster character of $E_x^{(n)}$ is independent of finite fields.

   \bdefinition
   (1) The $n$-th Chebyshev polynomials of the first kind is the polynomials $F_n(x)\in \Z[x]$ defined by 
      $$F_0(x)=1, F_1(x)=x,F_2(x)=x^2-2, F_{n+1}(x)=F_1(x)F_n(x)-F_{n-1}(x)\ for\ n\geq 2.$$
  
      (2) The $n$-th Chebyshev polynomials of the second kind is the polynomials $G_n(x)\in \Z[x]$ defined by 
      $$G_0(x)=1, G_1(x)=x,G_2(x)=x^2-1, G_{n+1}(x)=G_1(x)G_n(x)-G_{n-1}(x)\ for\ n\geq 2.$$
   \edefinition

   \blemma\label{lem5.3}
    We have
       $$F_n(X_{\delta})=X_{n\delta}\ \ \ and\ \ \ G_n(X_{\delta})=X_{E_x^{(n)}}.$$
    for some $x\in \P^1$ with degree 1.
    \bproof
    Note $X^{-(\begin{smallmatrix}r_1 \\ 0 \end{smallmatrix})}X^{-(\begin{smallmatrix}r_2 \\ 0 \end{smallmatrix})}=X^{-(\begin{smallmatrix}r_2 \\ 0 \end{smallmatrix})}X^{-(\begin{smallmatrix}r_1 \\ 0 \end{smallmatrix})}$ for $r_1,r_2\in \Z$. Consequently we can identity $X^{-(\begin{smallmatrix}n \\ 0 \end{smallmatrix})}$ with $z^{n}$ and $X^{-(\begin{smallmatrix}-n \\ 0 \end{smallmatrix})}$ with $z^{-n}$. By direct computations, we have
     $$F_2(X_{\delta})=(z+z^{-1})^2-2=z^2+z^{-2}=X_{2\delta}\ \ and \ \ G_2(X_{\delta})=z^2+1+z^{-2}=X_{E_x^{(2)}}.$$
    Assume $F_i(X_{\delta})=X_{i\delta}$ and  $G_i(X_{\delta})=X_{E_x^{(i)}}$  for $i\leq n$. Then
      \begin{align*}
      G_{n+1}(X_{\delta})&=X_{\delta}X_{E_x^{n}}-X_{E_x^{(n-1)}} \\
              &=(z+z^{-1})(z^n+z^{n-2}+\cdots +z^{-n})-(z^{n-1}+z^{n-3}+\cdots +z^{-(n-1)})\\
              &=z^{n+1}+z^{n-1}+\cdots t^{-n+1} +z^{-n-1}=X_{E_x^{(n+1)}}.
      \end{align*}
    $$F_{n+1}(X_{\delta})=X_{\delta}X_{n\delta}-X_{(n-1)\delta}=(z+z^{-1})(z^n+z^{-n})-(z^{n+1}+z^{-(n-1)})=X_{(n+1)\delta}.$$
    \eproof
   \elemma
  
   Next, we will show $X_{\delta}\in \sA(\varLambda,B)$, then as a result all $X_{n\delta}$ and $X_{E_x^{(n)}}$ will belong to $\sA(\varLambda,B)$. 
   
    \blemma\label{lem5.4}
      The following relations hold on $\sA(\varLambda,B)$:
       \begin{align}
        \label{id5.1} X_{\sO(2)}X_{\sO}&=\nu^2 X_{\sO(1)}^2+1.\\
        \label{id5.2} X_{\sO(3)}X_{\sO}&=\nu^2 X_{\sO(2)}X_{\sO(1)}+\nu^{-1} X_{\delta}.\\
        \label{id5.3} X_{\delta}X_{\sO}&=\nu X_{\sO(1)}+\nu^{-1}X_{\sO(-1)}.
       \end{align}
     \bproof
       The first one is the defining relation of $\sA(\varLambda,B)$.
      For any finite field $k$, the quantum cluster multiplication formula applied to $\sO(3)^k$ and $\sO^k$ is 
      $$(|k|^2-1)X_{\sO(3)^k}X_{\sO^k}=|k|^{\frac{3}{4}}(|k|^2-1)X_{\sO(2)^k\oplus \sO(1)^k}+|k|^{\frac{-1}{4}} \sum_{i=1}^{|k|+1}(|k|-1)X_{S_{x_i}^k}.$$
      where $S_{x_i}^k$ is the simple torsion sheaf supported on $x_i\in \P$ with $\mathrm{deg}(x_i)=1$. Note that $X_{S_{x_i}^k}=X_{\delta}$ for each $x_i$ and $X_{\sO(2)^k}X_{\sO(1)^k}=|k|^{\frac{1}{4}}X_{\sO(1)^k\oplus \sO(2)^k}$, we have
        $$X_{\sO(3)^k}X_{\sO^k}=|k|^{\frac{1}{2}}X_{\sO(2)^k}  X_{\sO(1)^k}+ |k|^{\frac{-1}{4}}X_{\delta}.$$
      Thus, in $\sA(\varLambda, B)$ we have
      $$X_{\sO(3)}X_{\sO}=\nu^2X_{\sO(2)}  X_{\sO(1)}+\nu^{-1} X_{\delta}.$$
      The proof of the third equation is similar due to Equation (\ref{eqn3.4}).
      \eproof
     \elemma

     Since $X_{\sO(3)}$ and $X_{\sO}$ belong to $\sA(\varLambda,B)$,  then $X_\delta=\nu X_{\sO(3)}X_{\sO}-\nu^3X_{\sO(2)}X_{\sO(1)}$ lies in $\sA(\varLambda,B)$. Applying the twisting operation $\sigma$ successfully to Identities (\ref{id5.1})  and (\ref{id5.2}), we have the following
    \begin{corollary}
      In $\sA(\varLambda,B)$, we have
      \begin{align}
           X_{\sO(l+2)}X_{\sO(l)}&=\nu^2 X_{\sO(l+1)}^2+1,\\
         X_{\sO(l+3)}X_{\sO(l)}&=\nu^2 X_{\sO(l+2)}X_{\sO(l+1)}+\nu^{-1} X_{\delta}.\\
       \label{eqn5.8}  X_{\delta}X_{\sO}&=\nu X_{\sO(1)}+\nu^{-1}X_{\sO(-1)}.
     \end{align}
    \end{corollary}

    \blemma\label{lem5.9}
    In $\sA^{tor}(\varLambda,B)$, we have
    \begin{align*}
      X_{n\delta}X_{m\delta}&=X_{m\delta}X_{n\delta}.\\
      X_{n\delta}X_{m\delta}&=X_{(n+m)\delta}+X_{(n-m)\delta}\ for\ n>m.\\
      X_{n\delta}X_{n\delta}&=X_{2n\delta}+2,\ \ for\ n\in \Z.
    \end{align*}
    Moreover $X_{\delta}^n$ is a $\Z$-linear combinations of $X_0$, $X_{\delta}$, $X_{2\delta}$, $\cdots$,  $X_{n\delta}$ with the coefficient of $X_{n\delta}$ being 1.
    \bproof
     The first and second statement is obtained by easy computation. We proceed an induction on $n$ to prove the third one. When $n=1$, the statement holds obviously. Assume it holds for $n-1$, then $X_{\delta}^n=X_{\delta}X_{\delta}^{n-1}$=$X_{\delta}(X_{(n-1)\delta}+a_1X_{(n-2)\delta}+\cdots +a_{n-2}X_{\delta}+a_{n-1})$ for some $a_i\in \Z$. Using the first statement, $X_{\delta}^n=X_{n\delta}+a_1X_{(n-1)\delta}+(a_1+1)X_{(n-2)\delta}+\cdots+b_n$ for $b_i\in \Z$.
     \eproof
     \elemma

   \bremark
   The Lemma \ref{lem5.9} above has been proved for the quantum cluster algebra $\sA(2,2)$ of Kronecker quiver in \cite[Proposition 6(1)]{Ding2012b}.
   \eremark

   \bproposition\label{prop5.10}
   Each one of the following sets forms a $\Z[\nu^{\pm1}]$-basis for $\sA^{tor}(\varLambda,B)$:
   $$\begin{aligned}
     \B_1^{tor}&=\{X_{r\delta}|r\in \N\},\\
     \B_2^{tor}&=\{X_{\delta}^r|r\in \N\},\\
     \B_3^{tor}&=\{X_{E_x^{(r)}}|r\in \N\}.\\
   \end{aligned}$$
   for some $x\in \P^1$ with degree 1.
   \bproof
   We have shown $\B_1^{tor}$ spans  $\sA^{tor}(\varLambda,B)$ in Lemma \ref{lem5.9}.  Since $X_{r\delta}=z^r+z^{-r}$ for distinct $r$ has different maximal degree, it follows that $\{X_{r\delta}|r\in \N\}$ is linearly independent. As a consequence, $\B_1^{tor}$ is a basis for $\sA(\varLambda,B)$.
   
   Note that $\B_1^{tor}$ and $\B_i^{tor}$ can be linearly represented by each other for $i=2,3$, which implies $\B_i^{tor}$ is also a basis.
   \eproof
  \eproposition
 
  \subsection{Bases of \texorpdfstring{$\sA(\varLambda,B)$}{Lx}}
  Since $\sA(\varLambda,B)$ is generated by $X_{\sO(l)}$ for $l\in \Z$, every element in $\sA(\varLambda,B)$ is a $\Z[\nu^{\pm}]$-linear combination of products of several $X_{\sO(l_i)}$. To find a basis for $\sA(\varLambda,B)$, we proceed by induction on the rank, i.e. the length (=s) of a product $\prod_{i=1}^s X_{\sO(l_i)}$.
  
  For the case when rank is 1, every element in $\sA(\varLambda,B)$ of rank 1 is a linear combination of $X_{\sO(l)}=X^{-(\begin{smallmatrix}l+1\\ -1\end{smallmatrix})}+\sum_{r\geq l}\nu^{-2(l+r)}[l+r+1]_{\nu^4}X^{-(\begin{smallmatrix}l+2r+1\\ 1\end{smallmatrix})}$, for $l\in \Z$.

  For the case when rank is 2,  we need to deal with $X_{\sO(n)}X_{\sO(m)}$ for any $n,m\in \Z$. Apply the operation of twisting, it suffices to deal with $X_{\sO}X_{\sO(n)}$ and  $X_{\sO(n)}X_{\sO}$ for $n\geq 0$. Denoted $z_n$ by $X_{n\delta}=F(X_{\delta})$,

  \bproposition \label{prop5.11}
   For $n\in \N$, $n\geq 0$, we have\ \ 
   \begin{align}
    \label{eqn5.9}&X_{\sO(2n)}X_{\sO}=\nu^{2n}X_{\sO(n)}^2+\sum_{l=0}^{n-1}\nu^{2(-n+2l+1)}\sum_{i=l+1}^{n}z_{2(n-i)}.\\
    \label{eqn5.10}&X_{\sO(2n+1)}X_{\sO}=\nu^{2n}X_{\sO(n+1)}X_{\sO(n)}+\sum_{l=0}^{n-1}\nu^{2(-n+2l)+1}\sum_{i=l+1}^{n}z_{2(n-i)+1}.
   \end{align}
   \bproof
   The proof is similar to the proof of \cite[Proposition 6(3)]{Ding2012b} using Lemma \ref{lem5.4}.
   \eproof
  \eproposition
  
  Applying the bar involution $\overline{\bullet}$ to Equations (\ref{eqn5.8}), (\ref{eqn5.9}) and (\ref{eqn5.10}), we have 
  \begin{align}
    X_{\sO}X_{\delta}&=\nu^{-1}X_{\sO(1)}+\nu X_{\sO(-1)}.\\
    X_{\sO}X_{\sO(2n)}&=\nu^{-2n}X_{\sO(n)}^2+\sum_{l=0}^{n-1}\nu^{-2(-n+2l+1)}\sum_{i=l+1}^{n}z_{2(n-i)}.\\
    \label{eqn5.13} X_{\sO}X_{\sO(2n+1)}&=\nu^{-2n}X_{\sO(n)}X_{\sO(n+1)}+\sum_{l=0}^{n-1}\nu^{-2(-n+2l)+1}\sum_{i=l+1}^{n}z_{2(n-i)+1}.
   \end{align}
  
   Define a subset $\C_r$ of $\sA(\varLambda,B)$ to be 
        $$\C_r=\{ X_{\sO(l)}^d X_{\sO(l+1)}^{r-d}|l\in \Z, 1\leq d\leq r\}.$$
   Set $\B^{vet}:=\bigcup_{r\geq 1} \C_r$.
 
   \btheorem\label{thm5.15}
    The set $\B^{vet}\bigcup \B_2^{tor}$ forms a $\Z[\nu^{\pm1}]$-basis for $\sA(\varLambda,B)$. 
   \bproof
    The set  linearly spans $\sA(\varLambda,B)$ following from Equations (\ref{eqn5.8})-(\ref{eqn5.13}). The linear independence of $\B^{vet}\bigcup \B_2^{tor}$ is induced from that each element in this set has distinct minimal degrees. Indeed, the minimal degrees of $X_{\sO(l)}$ are $-(\begin{smallmatrix}l+1\\ -1\end{smallmatrix})$ and $-(\begin{smallmatrix}3l+1\\ 1\end{smallmatrix})$. Thus the minimal degrees of $X_ {\sO(l)}^d X_{\sO(l+1)}^{r-d}$ are $-(r-d)(\begin{smallmatrix}l+1\\ -1\end{smallmatrix})-d(\begin{smallmatrix}l+2\\ -1\end{smallmatrix})$ and $-(r-d)(\begin{smallmatrix}3l+1\\ 1\end{smallmatrix})-d(\begin{smallmatrix}3l+4\\ 1\end{smallmatrix})$, which are $-(\begin{smallmatrix}rl+d+r\\ -r\end{smallmatrix})$ and $-(\begin{smallmatrix}3(lr+d)+r\\ r\end{smallmatrix})$ respectively.  When $r$ is fixed, two such elements share same minimal degree if and only if $rl+d=rl'+d'$. But $|d'-d|<r$ and $|l-l'|\geq 1$, it is impossible that there exist two different pairs $(l,d)$ and $(l',d')$ such that $rl+d=rl'+d'$. On the other hand the minimal degree of $X_\delta^n$ is $-(\begin{smallmatrix}-n\\ 0\end{smallmatrix})$. Hence different $(r,l,d)$ gives different minimal degrees ($X_\delta^n$ corresponding to $(0,0,n)$), it follows that the set $\{x\in \C_r|r\in \Z^+\}\cup \{X^n_{\delta}|n\in \N\}$ is linearly independent.
  \eproof
  \etheorem
  
\subsection{Bar-invariant bases}

  \begin{proposition}\ \ 

    \begin{itemize}
      \item [(i)] $\overline{\nu^{-1}X_{\sO(l)}X_{\sO(l+1)}}=\nu^{-1}X_{\sO(l)}X_{\sO(l+1)}$ and $\overline{X_{\sO(l)}X_{\sO(l)}}=X_{\sO(l)}X_{\sO(l)}$.
      \item [(ii)]$\nu^{d(r-d)}X_ {\sO(l)}^d X_{\sO(l+1)}^{r-d}$ is bar-invariant.
    \end{itemize} 
  \bproof
  The first statement is induced from Equation (\ref{id6.2}) and $X_{\sO(l)}$ is bar-invariant, i.e.,
  \begin{center}
    $\overline{X_{\sO(l)}X_{\sO(l+1)}}=\overline{X_{\sO(l+1)}}$ $\overline{X_{\sO(l)}}=\nu^{-2}X_{\sO(l)}X_{\sO(l+1)}$.
  \end{center} To show the second argument, by using $\overline{XY}=\overline{Y}$ $\overline{X}$, we have
    $$\overline{X_{\sO(l)}^d X_{\sO(l+1)}^{r-d}}=X_{\sO(l+1)}^{r-d} X_{\sO(l)}^d.$$
    Apply Equation \ref{id6.2} successfully on the left, we have
   $$\overline{\nu^{-d(r-d)}X_{\sO(l+1)}^{r-d} X_{\sO(l)}^d}=\nu^{-d(r-d)}X_{\sO(l)}^d X_{\sO(l+1)}^{r-d}.$$
   The proof is completed.
  \eproof
  \end{proposition}

  Set  $$\bar{\C}_r:=\{ \nu^{d(r-d)}X_{\sO(l)}^d X_{\sO(l+1)}^{r-d}|l\in \Z, 1\leq d\leq r\},$$
  and $\bar{\B}:=\bigcup_{r\geq 1} \bar{\C}_r$.  Combining Theorem \ref{thm5.15} with Proposition \ref{prop5.10}, we have
  \btheorem\label{thm5.14}
  Each one of the following sets gives rise to a  bar-invariant $\Z[\nu^{\pm1}]$-basis for $\sA(\varLambda,B)$:
  \begin{center}
    $ \B_1^{tor}\bigcup \bar{\B}^{vet}$, $\ \ $  $\B_2^{tor}\bigcup \bar{\B}^{vet}$, $\ \ $   $\B_3^{tor}\bigcup \bar{\B}^{vet}$.
  \end{center}
 \etheorem

 \bremark
  The isomorphism $\kappa:\sA(2,2)\to \sA(\varLambda,B)$ also preserves bar-invariant $\Z[\nu^{\pm1}]$-bases (see \cite[Corollary 9]{Ding2012b}).
 \eremark

  \appendix
  \section{Compatibility of exchange triangles}\label{ap}
  \subsection{Cluster categories of weighted projective lines}
   Let $k$ be any field. Recall $\X_{\bp,\bla}$ is the weighted projective line of $(\bp,\bla)$ given in Section 2.1.  Denoted by $\sA_k$ the hereditary category of coherent sheaves on $\X_{\bp,\bla}^k$.

   \bdefinition\label{def7.1}
    The canonical algebra $C(\bp,\bla)$ is defined to be the path algebra $kQ$ of the quiver $Q$ modulo following relations:
       $$\alpha_{i,1}\alpha_{i,2}\cdots \alpha_{i,p_i}=\alpha_{2,1}\cdots \alpha_{2,p_2}-\lambda_{i}\alpha_{1,1}\cdots \alpha_{1,p_1},$$
    for $i=3,\cdots, N$. Denote by $S$ the set consisting of the above relations.  Here $Q$ is given by
    $$
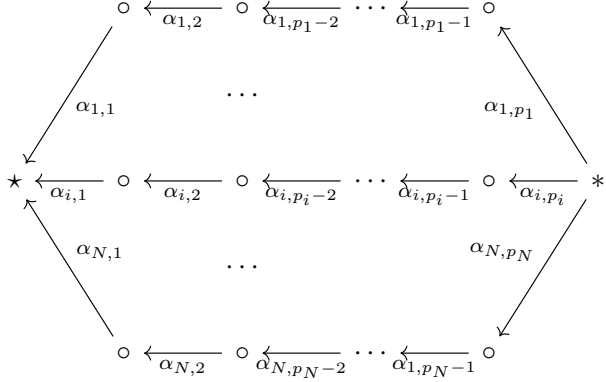
\begin{tikzcd}
      &\circ\arrow[ldd,"\alpha_{1,1}"]  &\circ\arrow[l,"\alpha_{1,2}"] &\cdots\arrow[l,"\alpha_{1,p_1-2}"] &\circ\arrow[l,"\alpha_{1,p_1-1}"] & \\
        & &\cdots & & &\\
     \star &\circ\arrow[l,"\alpha_{i,1}"] &\circ\arrow[l,"\alpha_{i,2}"] &\cdots\arrow[l,"\alpha_{i,p_i-2}"] &\circ\arrow[l,"\alpha_{i,p_i-1}"] &\ast\arrow[luu,"\alpha_{1,p_1}"]\arrow[l,"\alpha_{i,p_i}"]\arrow[ldd,"\alpha_{N,p_N}"']   \\
     & &\cdots  & & &\\
     &\circ\arrow[luu,"\alpha_{N,1}"']  &\circ\arrow[l,"\alpha_{N,2}"] &\cdots\arrow[l,"\alpha_{N,p_N-2}"] &\circ\arrow[l,"\alpha_{1,p_N-1}"]  & 
    \end{tikzcd}$$
   \edefinition
  
   Let $T^k=\bigoplus_{0\leq \vec{l}\leq \vec{c}} \sO(\vec{l})^k$ and $D^b(C(\bp,\bla))$ denoted by $D^b(\mathrm{mod}C(\bp,\bla))$. It is known that $T^k$ is a tilting object in $\sA_k$ and the derived functor $\R\Hom(T^k,-): D^b(\sA_k)\to D^b(C(\bp,\bla))$ is a derived equivalence, see \cite{Geigle1991,Chen2009}. Obviously, under this functor, the image of $\sO(j\vec{x_i})$  is the indecomposable projective $C(\bp,\bla)$-module $P_{i,j}$ for $1\leq i\leq N$ and $1\leq j\leq p_i-1$. $\R\Hom(T^k,\sO^k)=P_{\star}$ and $\R\Hom(T^k,\sO(\vec{c}))=P_{\ast}$.
   
   Let $\sC(\sA_k)$ be the cluster category $D^b(\sA_k)/\tau\circ [-1]$ of $\sA_k$. Since $\sA_k$ is derived equivalent to $\mathrm{mod}(C(\bp,\bla))$, the orbit category $D^b(C(\bp,\bla))/\tau\circ [-1]$ has a natural triangulated structure induced from $\sC(\sA_k)$ such that the projection functor $D^b(C(\bp,\bla))\to D^b(C(\bp,\bla))/\tau\circ [-1]$ is a triangle functor. We say $\sC(C(\bp,\bla)):=D^b(C(\bp,\bla))/\tau\circ [-1]$ is the cluster category of $C(\bp,\bla)$.

  \btheorem[{\cite[Theorem 6.8]{Buan2006a}}]\label{thm7.2}
  Let $T$ be a 2-Calabi-Yau triangulated category with a cluster-tilting object $T$. Let $T_i$ be indecomposable and $T = T_i \oplus \bar{T}$. Then there exists a unique indecomposable $T_i^*$ non-isomorphic to $T_i$ such that $\bar{T} \oplus T_i^*$ is cluster tilting. Moreover $T_i$ and $T_i^*$ are linked by the existence of exchange triangles
  $$T_i\stackrel{u}\lrw B \stackrel{v}\lrw T_i^* \lrw T_i[1]\qquad\text{and}\qquad 
   T_i^*\stackrel{u'}\lrw B' \stackrel{v'}\lrw T_i \lrw T^*_i[1]$$
  where $u$ and $u'$ are minimal left add $\bar{T}$-approximations and $v$ and $v'$ are minimal right add $\bar{T}$-approximations.
  \etheorem

  \btheorem[{\cite[Theorem 7.5]{Buan2006a}}]\label{thm7.3}
    Two indecomposable rigid objects $T_i$ and $T_i^*$ form an exchange pair if and only if 
        $$\dim_{\mathrm{End}_{\sC}T_i}\operatorname{Ext}^1_{\sC}(T_i,T_i^*)=1=\dim_{\mathrm{End}_{\sC}T_i^*}\operatorname{Ext}^1_{\sC}(T^*_i,T_i).$$
  \etheorem
  It is well known that for an indecomposable rigid object $\sF$ in $\tilde{\sA}_k$, we have that $\mathrm{End}_{\sA_k}(\sF)\cong k$ for any fields $k$. Moreover, since the global dimension of $\sA_k$ is $2$, we know that
      $$\frac{\mathrm{End}_{\sC_k}\sF}{\mathrm{rad}(\mathrm{End}_{\sC_k}\sF)}\cong \frac{\mathrm{End}_{\sA_k}\sF}{\mathrm{rad}(\mathrm{End}_{\sA_k}\sF).}$$ 
  Thus, $\mathrm{End}_{\sC_k}(\sF)\cong k$.

  \subsection{Base change functors}
    Let $q\in \N$ be prime, set $k=\F_q$ and $K=$. In the sequel, write  $A_{(r)}$ for $C(\bp,\bla)_{\F_{q^r}}$, $A_K$ for $C(\bp,\bla)_{K}$ and $\sA_{(r)}$ for $\sA_{\F_{q^r}}$. We have the following base change functors:
     $$\begin{tikzcd}
       \mathrm{mod}A_k\arrow[rr,"-\otimes_k \F_{q^r}"]\arrow[rrrr,bend left,"-\otimes_k K"] &  &\mathrm{mod}A_{(r)}\arrow[rr,"-\otimes_{\F_{q^r}} K"] & &\mathrm{mod}A_K.
      \end{tikzcd}$$
   Clearly, $-\otimes_k K=(-\otimes_{\F_{q^r}} K)\circ (-\otimes_k \F_{q^r})$. Since these base change functors are exact functors, we have following functors
       $$\begin{tikzcd}
       \sC(A_k)\arrow[rr,"-\otimes_k \F_{q^r}"]\arrow[rrrr,bend left,"-\otimes_k K"] &  &\sC(A_{(r)})\arrow[rr,"-\otimes_{\F_{q^r}} K"] & &\sC(A_k).
       \end{tikzcd}$$
   Set $M^{(r)}:=M\otimes_k \F_{q^r}$ and $M^{K}:=M\otimes_k K$ for $M\in \mathrm{mod}A_k$.

    \blemma\label{lem7.3}
     For $M, N\in \mathrm{mod}A_k$, then
     $$\Hom_{A_k}(M,N)\otimes_k K\cong \Hom_{A_K}(M^K,N^K).$$
     \bproof
     For any finite dimensional projective $A_k$-module $P$, we have an isomorphism
        $$\Hom_{A_k}(P,N)\otimes_k K\cong \Hom_{A_K}(P^K,N^K),$$
    since $A_K=A_k\otimes_k K$ and the above isomorphism holds for $A_k$, then for all direct summands of $A_k$. Indeed, $\Hom_{A_k}(A_k,N)\otimes_k K=N^K=\Hom_{A_K}(A_K,N^K)$.  Note that the global dimension of $A_k$ is 2, for $M\in A_k$, consider a projective resolution of $M$
        $$\begin{tikzcd}
       0\arrow[r] &P_2\arrow[r]  &P_1\arrow[r] &P_0\arrow[r] &M\arrow[r] &0.
        \end{tikzcd}$$
    Then  applying $-\otimes_k K$, we get a projective resolution of $M^K$.
    $$\begin{tikzcd}
      0\arrow[r] &P^K_2\arrow[r]  &P^K_1\arrow[r] &P^K_0\arrow[r] &M^K\arrow[r] &0.
       \end{tikzcd}$$
    Then for $N\in \mathrm{mod}A_k$, we have following commutative diagram with exact rows
       $$\begin{tikzcd}
        0\arrow[r] &\Hom_{A_k}(M,N)\otimes_k K\arrow[r]  &\Hom_{A_k}(P_0,N)\otimes_k K\arrow[r]\arrow[d] &\Hom_{A_k}(P_1,N)\otimes_k K\arrow[d]\\ 
        0\arrow[r] &\Hom_{A_K}(M^K,N^K)\arrow[r]  &\Hom_{A_K}(P^K_0,N^K)\arrow[r] &\Hom_{A_K}(P^K_1,N^K).\\ 
         \end{tikzcd}$$
      Here we use the fact $\Hom_{A_k}(M,N)$ is a $k$-linear vector space and $-\otimes_k K$ is exact. Note that the two vertical arrows are isomorphisms, it follows that 
         $$\Hom_{A_k}(M,N)\otimes_k K\cong \Hom_{A_K}(M^K,N^K).$$
     \eproof
    \elemma
    
    Let $K^b(\mathrm{proj}A_k)$ the bounded homotopy category of complexes of finite dimensional projective $A_k$-modules. Note that the global dimension of $C(\bp,\bla)_k$ is $2$, we known that $D^b(A_k)$ is triangle equivalent to $K^b(\mathrm{proj}A_k)$. For $P_{\bullet}\in K^b(\mathrm{proj}A_k)$, set $P_{\bullet}^K:=P_{\bullet}\otimes_k K$, which belongs to $K^b(\mathrm{proj}A_K)$. $\Hom_{A_k}(P_{\bullet},Q_{\bullet}):=\bigoplus_{i} \Hom_{A_k}(P_i,Q_i)$.
    \blemma\label{lem7.4}
     For $P_{\bullet}$, $Q_{\bullet}\in K^b(\mathrm{proj}A_k)$, we have an isomorphism
        $$\Hom_{K^b(\mathrm{proj}A_k)}(P_{\bullet},Q_{\bullet})\otimes_k K\cong \Hom_{K^b(\mathrm{proj}A_K)}(P_{\bullet}^K,Q_{\bullet}^K).$$
     \bproof
      Let 
         $$\Hom^{\bullet}(P_{\bullet},Q_{\bullet}):=\bigoplus_{i\in \Z}\Hom_{A_k}(P_{\bullet},Q_{\bullet}[i])$$
     be the complex of vector spaces with differential $d$ given by $d(f^i)=d_{Q}f^i-(-1)^if^id_{P}$ for $f^i\in \Hom_{A_k}(P_{\bullet},Q_{\bullet}[i])$. So we have that
         $$\Hom_{K^b(\mathrm{proj}A_k)}(P_{\bullet},Q_{\bullet})=H^0(\Hom^{\bullet}(P_{\bullet},Q_{\bullet})),$$
      and 
         $$\Hom_{K^b(\mathrm{proj}A_K)}(P^K_{\bullet},Q^K_{\bullet})=H^0(\Hom^{\bullet}(P^K_{\bullet},Q^K_{\bullet})).$$
      On the other hand, it can be checked that the following diagram is commutative:
        $$\begin{tikzcd}
        \Hom_{A_k}(P_{\bullet},Q_{\bullet}[-1])\otimes_k K\arrow[r]\arrow[d]  & \Hom_{A_k}(P_{\bullet},Q_{\bullet})\otimes_k K\arrow[r]\arrow[d] & \Hom_{A_k}(P_{\bullet},Q_{\bullet}[1])\otimes_k K\arrow[d]\\ 
        \Hom_{A_K}(P^K_{\bullet},Q^K_{\bullet}[-1])\arrow[r]  &\Hom_{A_K}(P^K_{\bullet},Q^K_{\bullet})\arrow[r] &\Hom_{A_K}(P^K_{\bullet},Q^K_{\bullet}[1]),
         \end{tikzcd}$$
      where vertical arrows are induced by isomorphisms
        $$\Hom_{A_k}(P_i,Q_j)\otimes_k K\cong \Hom_{A_K}(P^K_i,Q^K_j),$$
       as shown in Lemma \ref{lem7.3}. Hence we have that  
          $$H^0(\Hom^{\bullet}(P_{\bullet},Q_{\bullet})\otimes_k K)\cong \Hom_{K^b(\mathrm{proj}A_K)}(P^K_{\bullet},Q^K_{\bullet}).$$
      Notice that we deal with complexes of finite dimensional vector  spaces, it follows that 
        $$H^0(\Hom^{\bullet}(P_{\bullet},Q_{\bullet})\otimes_k K)\cong H^0(\Hom^{\bullet}(P_{\bullet},Q_{\bullet}))\otimes_k K.$$
     Thus, we get
       $$\Hom_{K^b(\mathrm{proj}A_k)}(P_{\bullet},Q_{\bullet})\otimes_k K\cong \Hom_{K^b(\mathrm{proj}A_K)}(P^K_{\bullet},Q^K_{\bullet}).$$
     \eproof
    \elemma

    Replacing $K$ with $\F_{q^r}$, we can show the above two statements over finite fields also hold.

    \blemma\label{lem7.5}\ \ 
    
    \begin{itemize}
      \item [(i)]For $M, N\in \mathrm{mod}A_k$, we have an isomorphism
        $$\Hom_{A_k}(M,N)\otimes_k \F_{q^r}\cong \Hom_{A_{(r)}}(M^{(r)},N^{(r)}).$$
      \item [(ii)] For $P_{\bullet}$, $Q_{\bullet}\in K^b(\mathrm{Proj}A_k)$, we have an isomorphism
      $$\Hom_{K^b(\mathrm{proj}A_k)}(P_{\bullet},Q_{\bullet})\otimes_k \F_{q^r}\cong \Hom_{K^b(\mathrm{proj}A_{(r)})}(P^{(r)}_{\bullet},Q^{(r)}_{\bullet}).$$
     \end{itemize}
    \elemma
   
    For an indecomposable object $\sF\in \sA_k$, there exists a complex $P_{\bullet}(\sF)$ of projective $A_k$-modules in  $K^b(\mathrm{proj}A_k)$ corresponding to $\sF$ under the derived functor $\R\Hom(T^k,-)$. Set $\sF^K\in \sA_K$ (resp. $\sF^{(r)}\in \sA_{\F_{q^r}}$) such that  $\R\Hom(T^K,P_{\bullet}(\sF)\otimes_k K)\cong \sF^K$ in $\sC(A_K)$ (resp.  $\R\Hom(T^{(r)},P_{\bullet}(\sF)\otimes_k \F_{q^r})\cong \sF^{(r)}$ in $\sC(A_{(r)})$). 

    \blemma\label{lem6.7}
    If $\sF\in \sA_{\F_q}$ is an indecomposable object, then both $\sF^{K}$ and $\sF^{(r)}$ are indecomposable in $\sA_K$ and $\sA_{\F_{q^r}}$ respectively, where $K=\bar{\F}_q$.
    \bproof
    Note that  $\sF$ is also indecomposable in $\sC(\sA_k)$ by \cite[Proposition 2.3]{Barot2010}, it follows that $\mathrm{End}_{\sC(\sA_k)}(\sF)$ is a local ring. Hence $\mathrm{End}_{\sC(\sA_k)}(\sF)\otimes_k K$ is also a local ring. Moreover by Lemma \ref{lem7.4} we have isomorphisms
     $$\mathrm{End}_{\sC(\sA_k)}(\sF)\otimes_k K\cong \mathrm{End}_{\sC(A_k)}(P_{\bullet}(\sF))\otimes_k K \cong \mathrm{End}_{\sC(A_k)}(P_{\bullet}(\sF)\otimes_k K). $$
     We can deduce that $\mathrm{End}_{\sC(\sA_K)}(\sF^K)\cong \mathrm{End}_{\sC(A_K)}(P_{\bullet}(\sF)\otimes_k K)$ is a local ring. One can show $\sF^{(r)}$ is an indecomposable object in the same way.
    \eproof
    \elemma

    Put everything together, by Theorems \ref{thm7.2} and \ref{thm7.3}, we obtain the following
    \btheorem\label{thm7.8}
     Let $(T_i,T_i^*)$ be an exchange pair in $\sC(\sA_k)$ with exchange triangles
      $$T_i\stackrel{u}\lrw B \stackrel{v}\lrw T_i^* \lrw T_i[1]\qquad\text{and}\qquad 
     T_i^*\stackrel{u'}\lrw B' \stackrel{v'}\lrw T_i \lrw T_i[1.]$$
     Then 
     \begin{itemize}
      \item [(i)]$(T_i^{(r)},T_i^{*(r)})$ is an exchange pair in $\sC(\sA_{(r)})$, whose exchange triangles are
      $$T_i^{(r)}\stackrel{u\otimes_k \F_{q^r}}\lrw B^{(r)} \stackrel{v\otimes_k \F_{q^r}}\lrw T_i^{*(r)} \lrw T_i^{(r)}[1],$$
   and
      $$ T_i^{*(r)}\stackrel{u'\otimes_k \F_{q^r}}\lrw B^{'(r)} \stackrel{v'\otimes_k \F_{q^r}}\lrw T_i^{(r)}\lrw T_i^{*(r)}[1].$$
      \item[(ii)]$(T_i^{K},T_i^{*K})$ is an exchange pair in $\sC(\sA_{K})$, whose exchange triangles are
      $$T_i^{K}\stackrel{u\otimes_k K}\lrw B^{K} \stackrel{v\otimes_k K}\lrw T_i^{*K} \lrw T_i^{K}[1],$$
   and
      $$ T_i^{*K}\stackrel{u'\otimes_k K}\lrw B^{'K} \stackrel{v'\otimes_k K}\lrw T_i^{K}\lrw T_i^{*K}[1].$$
     \end{itemize}
    \etheorem

    The cluster-tilting graph of $\sC(\sA_k)$ has as vertices the isomorphism classes of basic cluster-tilting objects of $\sC(\sA_k)$, while two vertices $T$ and $T'$ are connected by an edge if and only if they differ by precisely one indecomposable direct summand.

   \bcorollary\label{cor7.9}
   The cluster-tilting graph of $\sC(\sA_k)$ is connected if $k=\F_{q^r}$ or $\bar{\F}_q$, where $q$ is a prime and $r\geq 1$.
    \bproof
    By \cite[Theorem 1.2]{Fu2021a}, the cluster-tilting graph of $\sC(\sA_K)$ is connected for an algebraically closed field $K$. Since for any cluster-tilting object $T$ in $\sC(\sA_k)$, $T^K$ is a cluster-tilting object in $\sC(\sA_K)$ by Lemma \ref{lem7.4} and \ref{lem6.7}. It follows that the cluster-tilting graph of $\sC(\sA_K)$ is the same with the one of $\sC(\sA_k)$ from Theorem \ref{thm7.8}. The case when $k=\F_{q^r}$ is similar.
    \eproof
   \ecorollary

  \bibliographystyle{plain}
  \bibliography{myref}

\begin{thebibliography}{10}

\bibitem{Barot2010}
M.~Barot, D.~Kussin, and H.~Lenzing.
\newblock The cluster category of a canonical algebra.
\newblock {\em Transactions of the American Mathematical Society},
  362(8):4313--4330, 2010.

\bibitem{Berenstein2005}
Arkady Berenstein and Andrei. Zelevinsky.
\newblock Quantum cluster algebras.
\newblock {\em Advances in Mathematics}, 195(2):405--455, 2005.

\bibitem{Buan2011}
A.~B. Buan, O.~Iyama, I.~Reiten, and D.~Smith.
\newblock Mutation of cluster-tilting objects and potentials.
\newblock {\em American Journal of Mathematics}, 133(4):835--887, 2011.

\bibitem{Buan2006a}
Aslak~Bakke Buan, Robert Marsh, Markus Reineke, Idun Reiten, and Gordana
  Todorov.
\newblock Tilting theory and cluster combinatorics.
\newblock {\em Advances in Mathematics}, 204(2):572--618, 2006.

\bibitem{Buan2010}
Aslak~Bakke Buan, Robert~J. Marsh, and Dagfinn~F. Vatne.
\newblock Cluster structures from 2-{C}alabi-{Y}au categories with loops.
\newblock {\em Mathematische Zeitschrift}, 265(4):951--970, 2010.

\bibitem{Caldero2006a}
Philippe Caldero and Bernhard Keller.
\newblock From triangulated categories to cluster algebras. {II}.
\newblock {\em Annales Scientifiques de l'\'{E}cole Normale Sup\'{e}rieure.
  Quatri\`eme S\'{e}rie}, 39(6):983--1009, 2006.

\bibitem{Caldero2008}
Philippe Caldero and Bernhard Keller.
\newblock From triangulated categories to cluster algebras.
\newblock {\em Inventiones Mathematicae}, 172(1):169--211, 2008.

\bibitem{Chen2009}
Xiao-Wu Chen and Henning Krause.
\newblock Introduction to coherent sheaves on weighted projective lines.
\newblock 11 2009.

\bibitem{Chen2021}
Xueqing Chen, Ming Ding, and Haicheng Zhang.
\newblock The cluster multiplication theorem for acyclic quantum cluster
  algebras., 2021.

\bibitem{CrawleyBoevey1992}
William Crawley-Boevey.
\newblock Exceptional sequences of representations of quivers.
\newblock In {\em Proceedings of the {S}ixth {I}nternational {C}onference on
  {R}epresentations of {A}lgebras ({O}ttawa, {ON}, 1992)}, volume~14 of {\em
  Carleton-Ottawa Math. Lecture Note Ser.}, page~7. Carleton Univ., Ottawa, ON,
  1992.

\bibitem{Ding2013}
Ming Ding, Jie Xiao, and Fan Xu.
\newblock Integral bases of cluster algebras and representations of tame
  quivers.
\newblock {\em Algebras and Representation Theory}, 16(2):491--525, 2013.

\bibitem{Ding2012b}
Ming Ding and Fan Xu.
\newblock Bases of the quantum cluster algebra of the {K}ronecker quiver.
\newblock {\em Acta Mathematica Sinica (English Series)}, 28(6):1169--1178,
  2012.

\bibitem{Ding2012}
Ming Ding and Fan Xu.
\newblock Cluster characters for cyclic quivers.
\newblock {\em Frontiers of Mathematics in China}, 7(4):679--693, 2012.

\bibitem{Ding2020a}
Ming Ding, Fan Xu, and Haicheng Zhang.
\newblock Acyclic quantum cluster algebras via {H}all algebras of morphisms.
\newblock {\em Mathematische Zeitschrift}, 296(3-4):945--968, 2020.

\bibitem{SergeyFomin2002}
Sergey Fomin and Andrei Zelevinsky.
\newblock Cluster algebras i: Foundations.
\newblock {\em J. Amer. Math. Soc}, 15:497--529, 2002.

\bibitem{Fomin2003}
Sergey Fomin and Andrei Zelevinsky.
\newblock Cluster algebras ii: Finite type classification.
\newblock {\em Inventiones mathematicae}, 154(1):63--121, 2003.

\bibitem{Fu2021}
Changjian Fu and Shengfei Geng.
\newblock On cluster categories of weighted projective lines with at most three
  weights.
\newblock {\em Journal of Algebra}, 573:16--37, 2021.

\bibitem{Fu2021a}
Changjian Fu and Shengfei Geng.
\newblock On cluster-tilting graphs for hereditary categories.
\newblock {\em Advances in Mathematics}, 383:107670, 2021.

\bibitem{Fu2010}
Changjian Fu and Bernhard Keller.
\newblock On cluster algebras with coefficients and 2-{C}alabi-{Y}au
  categories.
\newblock {\em Transactions of the American Mathematical Society},
  362(2):859--895, 2010.

\bibitem{Fu2020}
Changjian Fu, Liangang Peng, and Haicheng Zhang.
\newblock Quantum cluster characters of hall algebras revisited.
\newblock {\em arXiv: Representation Theory}, 2020.

\bibitem{Geigle1987}
Werner Geigle and Helmut Lenzing.
\newblock A class of weighted projective curves arising in representation
  theory of finite-dimensional algebras.
\newblock In {\em Singularities, representation of algebras, and vector bundles
  ({L}ambrecht, 1985)}, volume 1273 of {\em Lecture Notes in Math.}, pages
  265--297. Springer, Berlin, 1987.

\bibitem{Geigle1991}
Werner Geigle and Helmut Lenzing.
\newblock Perpendicular categories with applications to representations and
  sheaves.
\newblock {\em Journal of Algebra}, 144(2):273--343, 1991.

\bibitem{Green1995}
James~A. Green.
\newblock Hall algebras, hereditary algebras and quantum groups.
\newblock {\em Inventiones Mathematicae}, 120(2):361--377, 1995.

\bibitem{Hubery2010}
Andrew Hubery.
\newblock Ringel-{H}all algebras of cyclic quivers.
\newblock {\em S\~{a}o Paulo Journal of Mathematical Sciences}, 4(3):351--398,
  2010.

\bibitem{Kapranov1997}
M.~M. Kapranov.
\newblock Eisenstein series and quantum affine algebras.
\newblock {\em Journal of Mathematical Sciences}, 84(5):1311--1360, 1997.

\bibitem{Kedzierski2013}
Dawid Kedzierski and Hagen Meltzer.
\newblock Schofield induction for sheaves on weighted projective lines.
\newblock {\em Communications in Algebra}, 41(6):2033--2039, 2013.

\bibitem{Keller2011}
Bernhard Keller.
\newblock Deformed {C}alabi-{Y}au completions.
\newblock {\em Journal f\"{u}r die Reine und Angewandte Mathematik. [Crelle's
  Journal]}, 654:125--180, 2011.
\newblock With an appendix by Michel Van den Bergh.

\bibitem{Lenzing2011}
Helmut Lenzing.
\newblock Weighted projective lines and applications.
\newblock In {\em Representations of algebras and related topics}, EMS Ser.
  Congr. Rep., pages 153--187. Eur. Math. Soc., Z\"{u}rich, 2011.

\bibitem{Lenzing2009}
Helmut Lenzing and Hagen Meltzer.
\newblock Exceptional pairs in hereditary categories.
\newblock {\em Communications in Algebra}, 37(8):2547--2556, 2009.

\bibitem{Meltzer1995}
Hagen Meltzer.
\newblock Exceptional sequences for canonical algebras.
\newblock {\em Archiv der Mathematik}, 64(4):304--312, 1995.

\bibitem{Palu2008}
Yann Palu.
\newblock Cluster characters for 2-{C}alabi-{Y}au triangulated categories.
\newblock {\em Universit\'{e} de Grenoble. Annales de l'Institut Fourier},
  58(6):2221--2248, 2008.

\bibitem{Palu2012}
Yann Palu.
\newblock Cluster characters {II}: a multiplication formula.
\newblock {\em Proceedings of the London Mathematical Society. Third Series},
  104(1):57--78, 2012.

\bibitem{Qin2012}
Fan Qin.
\newblock Quantum cluster variables via {S}erre polynomials.
\newblock {\em Journal f\"{u}r die Reine und Angewandte Mathematik. [Crelle's
  Journal]}, 668:149--190, 2012.
\newblock With an appendix by Bernhard Keller.

\bibitem{Rupel2011}
Dylan Rupel.
\newblock On a quantum analog of the {C}aldero-{C}hapoton formula.
\newblock {\em International Mathematics Research Notices. IMRN},
  14:3207--3236, 2011.

\bibitem{Schiffmann2004}
Olivier Schiffmann.
\newblock Noncommutative projective curves and quantum loop algebras.
\newblock {\em Duke Mathematical Journal}, 121(1):113--168, 2004.

\bibitem{Schiffmann2012a}
Olivier Schiffmann.
\newblock Lectures on {H}all algebras.
\newblock In {\em Geometric methods in representation theory. {II}}, volume~24
  of {\em S\'{e}min. Congr.}, pages 1--141. Soc. Math. France, Paris, 2012.

\bibitem{Xiao2010}
Jie Xiao and Fan Xu.
\newblock Green's formula with {$\mathbb{C}^*$}-action and {C}aldero-{K}eller's
  formula for cluster algebras.
\newblock In {\em Representation theory of algebraic groups and quantum
  groups}, volume 284 of {\em Progr. Math.}, pages 313--348.
  Birkh\"{a}user/Springer, New York, 2010.

\bibitem{Xu2010}
Fan Xu.
\newblock On the cluster multiplication theorem for acyclic cluster algebras.
\newblock {\em Transactions of the American Mathematical Society},
  362(2):753--776, 2010.

\bibitem{Zhang1996}
Pu~Zhang.
\newblock Triangular decomposition of the composition algebra of the
  {K}ronecker algebra.
\newblock {\em Journal of Algebra}, 184(1):159--174, 1996.

\end{thebibliography}

 \end{document}